\def\draft#1{}                      

\documentclass[11pt]{article}

\newtheorem{theorem}{Theorem}[section]
\newtheorem{lemma}[theorem]{Lemma}
\newtheorem{prop}[theorem]{Proposition}
\newtheorem{example}[theorem]{Example}
\newtheorem{remark}[theorem]{Remark}

\long\def\rests#1{}

\def\noi{\noindent}

\def\pf{\noi{\bf Proof.\ \,}}
\def\eop{\hfill\framebox[2.4mm][t1]{\phantom{x}} \vskip 0.15cm } 

\def\sss#1{\if#1..\ \else\if#1,,\ \else\ #1\fi\fi} 

\def\voa#1{vertex operator algebra\sss{#1}}
\def\voas#1{vertex operator algebras\sss{#1}}
\def\svoa#1{vertex operator super algebra\sss{#1}}
\def\svoas#1{vertex operator super algebras\sss{#1}}
\def\subvoa#1{vertex operator subalgebra\sss{#1}}
\def\subsvoa#1{vertex operator super subalgebra\sss{#1}}
\def\subvoas#1{vertex operator subalgebras\sss{#1}}

\def\VF{V_{\rm Fermi}}
\def\Aut{{\rm Aut}}

\def\C{{\bf C}}
\def\R{{\bf R}}
\def\Z{{\bf Z}}

\def\tr{{\rm tr}}

\def\h{{\bf h}}
\def\k{{\bf k}}
\def\s{{\bf s}}
\def\o{\,{\rm o}}

\title{Conformal Designs based on Vertex Operator Algebras}

\author{Gerald~H\"ohn\thanks{Department of Mathematics, 
Kansas State University, 138 Cardwell Hall, Manhattan, KS 66506-2602, USA.
 E-mail: {\tt gerald@math.ksu.edu} \newline
Supported by Kansas NSF EPSCoR grant NSF32239KAN32240.
}}

\date{}

\begin{document}

\maketitle

\abstract{We introduce the notion of a conformal design based on a 
vertex operator algebra. This notation is a natural 
analog of the notion of block designs or spherical designs when the
elements of the design are based on self-orthogonal binary codes or 
integral lattices, respectively. It is shown that the subspaces of fixed degree 
of an extremal self-dual vertex operator algebra form conformal $11$-, $7$-, or $3$-designs,
generalizing similar results of Assmus-Mattson and Venkov for extremal
doubly-even codes and extremal even lattices. Other examples are coming from
group actions on vertex operator algebras, the case studied first by Matsuo.
The classification of conformal $6$- and $8$-designs is investigated. Again, our results 
are analogous to similar results for codes and lattices. }


\section{Introduction}

In the past, it has been a fruitful approach  to generalize concepts
known for codes and lattices to vertex operator algebras and then to show
that analogous results hold in this context. Important examples are
the construction of the Moonshine module~\cite{Bo-ur,FLM} or the modular
invariance of the genus-$1$ correlation functions~\cite{Zhu-dr}. Other examples
involving \svoas are given in~\cite{Ho-dr,Ho-shadow}. A comprehensive analysis of
these analogies can be found in~\cite{Ho-kleinian}, where also a fourth step
in this analogy, codes over the Kleinian four-group ${\bf Z}_2\times {\bf Z}_2$, was introduced. 

Two well-studied combinatorial structures are block designs and spherical designs
and many examples of such designs are related to codes and lattices. A notion
analogous to  block designs and spherical designs in the context of \voas however
has been missing so far. In this paper, we introduce the notation of a {\it conformal 
design based on a vertex operator algebra\/} and prove several results analogous
to known results for block and spherical designs. 

\medskip

The paper is organized as follows. In the rest of the introduction we discuss
the concepts of  block and spherical $t$-designs to motivate our definition
of conformal $t$-designs. Details of the definition will be discussed in the next 
section. Section~2 also contains several basic results.
In particular, we prove that a vertex operator algebra $V$ with a compact automorphism
group $G$ leads to conformal designs if the fixed point \voa $V^G$ has large minimal
weight (Theorem~\ref{groupdesigns}).
The minimal weight of a \svoa $V$ is the smallest $m>0$ with $V_m\not=0$.
We also describe the construction of derived designs (Theorem~\ref{derivedparts}). 
In Section~3, we prove that the homogeneous parts
of an extremal \voa are conformal $11$-, $7$-, or $3$-designs depending on the residue 
class  $\bmod \ 24$ of the central charge of $V$.
Extremal \voas have been introduced
in~\cite{Ho-dr} and are \voas $V$ with only one irreducible module up to isomorphism
such that the first coefficients of the graded trace of $V$ are as small as possible. 
The Moonshine module is an example of an extremal \voa supporting $11$-designs. 
This theorem is the analog of theorems of Assmus-Mattson~\cite{AsMa} (see also~\cite{Ba-harm}) and 
Venkov~\cite{VE-design} for extremal doubly-even codes and extremal even lattices, 
respectively. In the final section, we study the classification of conformal 
$6$- and $8$-designs $V_m$ supported by \svoas $V$ of minimal weight $m\leq 2$. 

Our result for conformal $6$-designs (Theorem~\ref{classfic7}) is analogous to 
(but a little weaker than) similar results for 
codes~\cite{LaLa-codesign} and lattices~\cite{Ma-ladesign}: For $m=1/2$ and $m=1$,
the only examples are the single Fermion \svoa $V_{\rm Fermi}$ and the two lattice
\voas $V_{A_1}$ and $V_{E_8}$, respectively. For $m=3/2$, the allowed central charges are $c=16$
and $c=23\frac{1}{2}$. 
In the later case the shorter Moonshine module $V\!B^{\natural}$~\cite{Ho-dr} is a known example.
For $m=2$, the allowed central charges are $c=8$, $16$,
$23\frac{3}{35}$, $23\frac{1}{2}$, $24$,  $32$,  $32\frac{4}{5}$, $33\frac{5}{7}$, $34\frac{4}{7}$,
$40$, $1496$ and examples are known for $c=8$, $16$, $23\frac{1}{2}$, $24$ and $32$.
For conformal $8$-designs, one obtains (Theorem~\ref{classfic8}) $m=2$ and $c=24$
with the Moonshine module $V^{\natural}$ the only known example.

\smallskip
I like to thank J.-C.~Puchta and C.~Schweigert for answering questions and H.~Yamauchi for
his suggestions and his carefully reading of an early version of this paper.   

\bigskip

A {\it block design\/} $X$ of type $t$-$(n,k,\lambda)$ (or $t$-design, for short) 
is a set of size-$k$-subsets of 
a set $M$ of size $n$ such that each size-$t$-subset of $M$ is contained in the
same number $\lambda$ of sets from $X$. A subset $K$ of $M$ can be 
identified with the vector $(x_{q})_{q\in M}\in {\bf F}_2^M$, where $x_q=1$ if $q\in K$ and
$x_q=0$ otherwise. 
Let ${\rm Hom}(m)$ be the complex vector space 
of homogeneous  {\it polynomial functions\/} of degree $m$ on ${\bf F}_2^M$ in the variables
$t_p$, $p\in M$, where $t_p((x_q)_{q\in M})=x_p$ for $(x_q)_{q\in M} \in {\bf F}_2^M$.
A basis of ${\rm Hom}(m)$ is formed by the  $n \choose m$  monomials
$\prod_{p\in N} t_p$ where $N$ is a size-$m$-subset of $M$. 

There is a natural action of ${\Aut M}\cong  {\rm S}_n$, the symmetric group of degree~$n$, 
on $M$, ${\bf F}_2^M$, and ${\rm Hom}(m)$. The condition for a set $X$ of size-$k$-subsets of $M$
to be a $t$-design is now equivalent to (one may choose for $f$ the monomial
$\prod_{p\in N} t_p$ where $N \subset M$ with $|N|=s$)  that for all $0\leq s\leq t$ and  
$f\in  {\rm Hom}(s)$ and for all $g\in {\rm S}_n$ the following equation holds
\begin{equation}\label{b-designinv}
\sum_{B\in X}f(B)=\sum_{B\in X}f(g(B)).
\end{equation}

The vector spaces ${\rm Hom}(m)$ can be decomposed into a direct sum of irreducible components under the action
of ${\rm S}_n$. Let $\pi$ be the projection of ${\rm Hom}(m)$ onto the trivial
component. By averaging equation~(\ref{b-designinv}) over all $g\in {\rm S}_n$ one finds that
the above condition on $X$ is equivalent to that 
for all $0\leq s\leq t$ and  $f\in  {\rm Hom}(s)$ one has
\begin{equation}\label{b-designpro}
\sum_{B\in X}f(B)=\sum_{B\in X}\pi(f)(B).
\end{equation}

The decomposition of ${\rm Hom}(m)$ into irreducible components can be done explicitly:
The kernel of $\Delta:{\rm Hom}(m)\longrightarrow  {\rm Hom}(m-1)$, 
$f\mapsto \sum_{p\in M} \frac{ \partial f}{\partial t_p}$, forms the irreducible constituent 
${\rm Harm}(m)$ of
``discrete'' harmonic polynomials of degree~$m$. One has the decomposition ${\rm Hom}(m)=\bigoplus_{i=0}^{m}
{\rm Harm}'(i)$, where ${\rm Harm}'(i)$ is the image of ${\rm Harm}(i)$ in  ${\rm Hom}(m)$
under the multiplication with $(\sum_{p\in M}t_p)^{m-i}$. This gives the following
characterization of $t$-designs:  A set $X$ of size-$k$-subsets of $M$ is a $t$-design
if and only if for all $1\leq s\leq t$ and $f\in {\rm Harm}(s)$ one has
\begin{equation}\label{b-designharm}
\sum_{B\in X}f(B)=0.
\end{equation}

\medskip

A {\it spherical $t$-design\/} (cf.~\cite{DGS,GoSe}) is a finite subset $X$ of 
a sphere $S(r)$ of radius~$r$ 
around $0$ in ${\bf R}^n$ (usually the unit sphere) such that $X$ can be used
to integrate all polynomials $f$ of degree $s\leq t$ on $S(r)$ exactly by averaging
their values on $X$, i.e., the following equation holds:
\begin{equation}\label{s-designdef}
\frac{1}{|X|}\sum_{B\in X}f(B)=\frac{1}{{\rm Vol}(S(r))}\int_{S(r)}f\ \omega,
\end{equation}
where $\omega$ is the canonical volume form on the Riemannian manifold $S(r)$.
 
Let ${\rm Hom}(m)$ be the complex vector space 
of  homogeneous  polynomials of degree $m$  in the variables
$x_i$, where the $x_i$, $i=1$, $\ldots$, $n$, are  identified with
the standard coordinate functions on  ${\bf R}^n$. 
One has $\dim {\rm Hom}(m)={n+m-1\choose n-1}$.
There is a natural action of the orthogonal group ${\rm O}(n)$ 
on ${\bf R}^n$ and ${\rm Hom}(m)$. The condition for a set $X\subset S(r)$
to be a spherical $t$-design is now equivalent to (the right hand side of
equation~(\ref{s-designdef}) is obviously ${\rm O}(n)$-invariant)
that for all $0\leq s\leq t$ and  
$f\in  {\rm Hom}(s)$ and for all $g\in {\rm O}(n)$ the following equation holds
\begin{equation}\label{s-designinv}
\sum_{B\in X}f(B)=\sum_{B\in X}f(g(B)).
\end{equation}

The vector spaces ${\rm Hom}(m)$ can again be decomposed into a direct sum of 
irreducible components under the action
of ${\rm O}(n)$. Let $\pi$  be the projection of ${\rm Hom}(m)$ onto the trivial
component. By averaging equation~(\ref{s-designinv}) over all $g\in {\rm O}(n)$
using an invariant measure, one finds that
the above condition on $X$ is equivalent to that 
for all $0\leq s\leq t$ and  $f\in  {\rm Hom}(s)$ one has
\begin{equation}\label{s-designpro}
\sum_{B\in X}f(B)=\sum_{B\in X}\pi(f)(B).
\end{equation}

The decomposition of ${\rm Hom}(m)$ into irreducible components can also be done explicitly:
The kernel of $\Delta:{\rm Hom}(m)\longrightarrow  {\rm Hom}(m-2)$, 
$f\mapsto \sum_{l=1}^n \frac{ \partial^2f}{\partial x^2_l}$, forms the irreducible constituent 
${\rm Harm}(m)$ of
harmonic polynomials of degree~$m$. One has the decomposition ${\rm Hom}(m)=\bigoplus_{i=0}^{[m/2]}
{\rm Harm}'(m-2i)$, where ${\rm Harm}'(m-2i)$ is the image of ${\rm Harm}(m-2i)$ in  ${\rm Hom}(m)$
under the multiplication with $(\sum_{l=1}^n x^2_l)^{i}$. This gives the following
characterization of spherical $t$-designs:  A finite subset $X\subset S(r)$  is a spherical
$t$-design if and only if for all $1\leq s\leq t$ and $f\in {\rm Harm}(s)$ one has
\begin{equation}\label{s-designharm}
\sum_{B\in X}f(B)=0.
\end{equation}

\medskip

In the definition of a {\it conformal $t$-design based on a \voa $V$\/} we will replace the groups 
${\rm S}_n$ and ${\rm O}(n)$ by the Virasoro algebra of central charge~$c$. 
The r\^ole of  ${\rm Hom}(m)$ will be the degree-$m$-part $V_m$ of a \voa $V$. 
For $X$, we will take a homogeneous part of a module of $V$ and the evaluation 
$f[X]=\sum_{B\in X}f(X)$ will be replaced by
$$f[X]=\tr|_{X}\o(f),$$
where $f\in V_m$ and $\o(f)$ is the coefficient $v_{(k)}$ of the vertex operator 
$Y(v,z)=\sum_{n\in {\bf Z}} v_{(n)}z^{-n-1}$ which maps $X$ into itself.
The definition of a conformal $t$-design is then the same as in (\ref{b-designpro})
and (\ref{s-designpro}) or, equivalently, as in (\ref{b-designharm})
and (\ref{s-designharm}).


\section{Definition and basic properties}\label{basicdef}

The Virasoro algebra is the complex Lie algebra spanned by $L_n$, $n\in {\bf Z}$, 
and the central element $C$ with Lie bracket
\begin{equation}\label{virasoro} 
[L_m,L_n]=(m-n)L_{m+n}+\frac{m^3-m}{12}\delta_{m+n,0}\,C
\end{equation}
where $\delta_{k,0}=1$ if $k=0$ and $\delta_{k,0}=0$ otherwise.

For a pair $(c,h)$ of real numbers the Verma module $M(c,h)$ is 
a representation of the Virasoro algebra generated by 
a highest weight vector $v\in M(c,h)$ with $Cv=c$, $L_0v=hv$ and $L_nv=0$ for $n\geq 1$.
The number $c$ is called the central charge and $h$ is called the
conformal weight of the module. The set
$$\{L_{-m_1}\cdots L_{-m_k}v\mid k\,,m_1,\,\ldots,\,m_k\in \Z_{\geq 0},\, m_1\geq\cdots\geq m_k\geq 1\}$$
forms a basis of $M(c,h)$. The vector space $M(c,h)$ is graded by the eigenvalues of $L_0$ and
the vector $L_{-m_1}\cdots L_{-m_k}v$ is homogeneous of degree $m_1+\cdots+m_k+h$.

For $h=0$, the module $M(c,0)$ has a quotient isomorphic to $M(c,0)/M(c,1)$ with a basis
which can be identified with 
$$\{L_{-m_1}\cdots L_{-m_k}v\mid k,\,m_1,\,\ldots,\,m_k\in \Z_{\geq 0},\, m_1\geq\cdots\geq m_k\geq 2\}.$$

A vector $w$ of $M(c,h)$ or $M(c,0)/M(c,1)$ is said to be a singular vector if $L_mw=0$ for
all $m\geq 1$. The Kac-determinant formula shows that a
module $M(c,h)$ for $h\not=0$ contains a singular vector of degree up to $n$ 
if $h_{p,q}(c)=h$ for $pq\leq n$ where 
$$h_{p,q}(c)=\frac{((m+1)p-mq)^2-1}{4m(m+1)}\ \ \hbox{with}\ \ m=-\frac{1}{2}\pm\frac{1}{2}\sqrt{\frac{25-c}{1-c}}.$$

The modules $M(c,0)/M(c,1)$ contain a singular vector
of degree up to $n$ if and only if the central charge $c$ is a zero of the following 
normalized polynomial $D_n(c)$:
\begin{eqnarray}\label{kacdet} \nonumber
D_2(c)& = & c \\ \nonumber
D_4(c)& = &  c(5c+22) \\ \nonumber
D_6(c) & = & c(2c-1)(5c+22)(7c+68) \\ 
D_8(c) & = & c(2c-1)(3c+46)(5c+3)(5c+22)(7c+68).
\end{eqnarray}

\medskip

A vertex operator algebra $V$ over the field of complex numbers is a  complex vector space 
equipped with a linear map $Y:V\longrightarrow {\rm End}(V)[[z,z^{-1}]]$ 
and two nonzero vectors ${\bf 1}$ and $\omega$ in $V$ satisfying certain axioms; 
cf.~\cite{Bo-ur,FLM}. For $v\in V$ one writes
$$Y(v,z)=\sum_{n\in\Z}a_{(n)}z^{-n-1}.$$

The {\it vacuum vector\/} vector ${\bf 1}$ satisfies ${\bf 1}_{(-1)}={\rm id}_V$ and ${\bf 1}_{(n)}=0$
for $n\not=-1$. The coefficients $L_n=\omega_{(n+1)}$ for the {\it Virasoro vector\/} $\omega$
are satisfying the Virasoro relation~(\ref{virasoro}) with $C=c\cdot{\rm id}_V$ for a complex number $c$
called the {\it central charge\/} of $V$.

The operator $L_0$ give rise to a grading $V=\bigoplus_{n\in \Z} V_n$, where $V_n$ denotes
the eigenspace of $L_0$ with eigenvalue $n$, called the degree. $V_n$ is supposed to be finite-dimensional
and we assume $V_n=0$ for $n<0$ and $V_0=\C\cdot{\bf 1}$. For $v\in V_n$ the operator $v_{(n-1)}$ is
homogeneous of degree~$0$. We define $\o(v)=v_{(n-1)}$.

We will make use several times of the {\it associativity relation\/}
$$(u_{(m)}v)_{(n)}=\sum_{i\geq 0} (-1)^i{m \choose i }\bigl(
u_{(m-i)}v_{(n+i)}-(-1)^m v_{(m+n-i)}u_{(i)}\bigr)$$
for elements $u$, $v\in V$.

\smallskip
For the notion of admissible and (ordinary) module we refer to~\cite{DLM-twistedrep}.
A \voa is called {\it rational\/} if every admissible module is completely reducible. In this case
there are only finitely many irreducible admissible modules up to isomorphism and
every irreducible admissible module is an ordinary module. A \voa is called {\it simple\/} 
if it is irreducible as a module over itself. 

For an irreducible module $W$ there exists an $h$ such that $W=\bigoplus_{n\in \Z_{\geq 0}} W_{n+h}$,
where  the degree $n$ subspace $W_n$ is  again  the eigenspace of $L_0$ of eigenvalue $n$, with 
$W_h\not= 0$. We call  $h$ the {\it conformal weight\/} of the module $W$.

\medskip

The graded trace of an element $v\in V_n$ on a module $W$ of conformal weight~$h$ 
is defined by 
$$\chi_W(v,q)=q^{-c/24}\sum_{k\in \Z_{\geq 0}}\tr|_{W_{k+h}}\o(v)q^{k+h}.$$
For $v={\bf 1}$, we call $\chi_W(q)=\chi_W({\bf 1},q)$ the {\it character\/} of $W$.
If $V$ is assumed to be rational and satisfying the {\it $C_2$-cofiniteness condition\/}
$V/\hbox{Span}\{x_{(-2)}y\mid x,y\in V\}<\infty$ it is a result of Zhu that
$\chi_W(v,q)$ is a holomorphic function on the complex upper half plane
in the variable $\tau$ for $q=e^{2\pi i \tau}$. We assume  in this paper that the
$C_2$-cofiniteness condition is satisfied.

Furthermore, for $v \in V_n$ a highest weight vector for the Virasoro algebra,
the family $\{\chi_W(v,q)\}_W$, where $W$ runs through the isomorphism classes of
irreducible $V$-modules $W$, transforms as a vector-valued modular form of weight~$n$ 
for the modular group ${\rm PSL}_2(\Z)$ acting on the upper half plane in the usual
way.

Given a \voa $V$, Zhu defined a new \voa on the same underlying vector space 
with vertex operator $Y[v,z]=Y(v,e^z-1)e^{z{\rm wt}(v)}$ for homogeneous 
elements $v\in V$. The vacuum element is the same as the original one and
the new Virasoro element is $\tilde\omega=\omega-\frac{c}{24}$. We let
$Y[\tilde \omega,z]=\sum_{n\in \Z}L_{[n]}z^{-n-2}$. The new Virasoro
algebra generator $L_{[0]}$ introduces a new grading $V=\bigoplus_{n=0}^{\infty} V_{[n]}$
on~$V$ and similar for modules $W$. 
One has $\bigoplus_{n\leq N}W_n=\bigoplus_{n\leq N} W_{[n]}$ for all $N\in \R$.
In particular, the Virasoro highest weight vectors for both Virasoro algebras
are the same.

It was shown by H.~Li that a \voa $V$ as above has a unique normalized invariant bilinear form
$(\,.\,,\,.\,)$ given by
$$(u,v){\bf 1}={\rm Res}_{z=0}(z^{-1}Y(e^{zL_1}(-z^{-2})^{L(0)}u,z^{-1})v$$
for elements $u$, $v$ in $V$ with normalization $({\bf 1},{\bf 1})=1$ provided that $L_1V_1=0$.  

\medskip

We {\it assume\/} that the \voas $V$ in this paper are
isomorphic to a direct sum of highest weight modules for the Virasoro algebra,
i.e., one has $$V=\bigoplus_{i\in I } M_i,$$
where each $M_i$ is a quotient of a Verma modules $M(c,h)$ with $h\in \Z_{\geq 0}$.
One has therefore a natural decomposition   
\begin{equation}\label{natdecomp}
V= \bigoplus_{h=0}^{\infty}\overline{M}(h)
\end{equation}
where $\overline{M}(h)$ is a direct sum of finitely many quotients of the Verma module $M(c,h)$.
The module $\overline{M}(0)$ is the \subvoa of $V$ generated by $\omega$ which we denote
also by $V_{\omega}$ and is therefore a quotient of $M(c,0)/M(c,1)$.
The smallest $h>0$ for which $\overline{M}(h)\not= 0$ is called
the {\it minimal weight\/} of $V$ and denoted by $\mu(V)$. (If no such $h>0$ exists, we let
$\mu(V)=\infty$.) We note that our assumption implies that $L_1V_1=0$.

In particular, the decomposition~(\ref{natdecomp}) gives us the natural projection map
$$\pi: V \longrightarrow V_{\omega}$$
with kernel $\bigoplus_{h>0}^{\infty}\overline{M}(h)$.

The decomposition $V= \bigoplus_{h=0}^{\infty}\overline{M}(h)$ is the same for the
modified  Virasoro algebra with generators $L_{[n]}$.

\medskip

In Section~\ref{classresult}, we study \svoas. For the full definition we refer to~\cite{Ho-dr,Kac-VOA}.
We note that a \svoa is a super vector space $V=V_{(0)}\oplus V_{(1)}$
where the even part $V_{(0)}=\bigoplus_{n=0}^{\infty} V_n$ is a \voa and 
the odd part $V_{(1)}=\bigoplus_{n=0}^{\infty} V_{n+1/2}$ is a $V_{(0)}$-module 
with a conformal weight $h\in \Z_{\geq 0}+\frac{1}{2}$. The minimal weight $\mu(V)$ of 
$V$ is defined to be the minimum of $\mu(V_{(0)})$ and the conformal weight $h$ of $V_{(1)}$.
A \svoa $V$ is called rational in this paper if the even \subvoa $V_{(0)}$ is rational. 

\bigskip

As explained in the introduction, the following definition is motivated by
analogous definitions of block designs and spherical designs.

\medskip

\noi{\bf Definition:}  Let $V$ be \voa of central charge $c$ and let $X$ be 
a degree $h$ subspace of a module of $V$. For a positive integer $t$ one calls
$X$ a {\it conformal design of type $t$-$(c,h)$\/}
or {\it conformal $t$-design,\/} for short,
if for all $v\in V_n$ where $0\leq n\leq t$ one has
$$ \tr|_{X}\o(v)= \tr|_{X}\o(\pi(v)).$$

\smallskip

The following two observations are clear:
\begin{remark} If $X$ is a conformal $t$-design based on $V$, it is also a conformal $t$-design
based on an arbitrary \subvoa of $V$.\eop
\end{remark}

\begin{remark}
A conformal $t$-design is also a conformal $s$-design for all
integers \hbox{$1\leq s< t$.} \phantom{xx} \eop
\end{remark}

\begin{theorem}\label{equivalence}
Let $X$ be the homogeneous subspace of a module of a \voa $V$.
The following conditions are equivalent:
\begin{itemize}
\item[(i)] $X$ is a conformal $t$-design.
\item[(ii)] For all homogeneous  $v\in \ker \pi =\bigoplus_{h>0} \overline{ M}(h)$ of degree $n\leq t$, one
has  $\tr|_{X}\o(v)=0$. 
\end{itemize}
\end{theorem}
\pf For $v\in V$ one has $v-\pi(v) \in  \bigoplus_{h>0} \overline{ M}(h)$.  \eop  

If $V$ has minimal weight $\mu$, then the homogeneous subspaces of any module of $V$ are
conformal $t$-designs for all $t=1$, $2$, $\ldots$, $\mu-1$. We call
such a conformal $t$-design {\it trivial\/}. In particular,
if $V$ is isomorphic to $V_{\omega}$, i.e., has minimal weight infinity,
all conformal designs based on $V$ are trivial.

\begin{theorem}\label{highestweight}
Let $V$ be a  \voa and let $N$ be a $V$-module graded by $\Z +h$. 
The following conditions are equivalent:
\begin{itemize}
\item[(i)] The homogeneous subspaces
$N_n$ of $N$ are conformal $t$-designs based on $V$ for $n\leq h$.
\item[(ii)] For all Virasoro highest weight vectors $v\in V_s$ with $0<s\leq t$
and all $n\leq h$ one has 
$$\tr|_{N_n}\o(v)=0.$$
\end{itemize}
\end{theorem}
\pf It is clear that (i) implies (ii). 
Assume now that (ii) holds. Let $v\in V$ be a Virasoro highest weight vector and
let $w=L_{[-i_1]}\cdots L_{[-i_k]}v$ with positive integers $i_1$, $\ldots$, $i_k$
be in the Virasoro highest weight module generated by $v$. 
In~\cite{Zhu-dr}, Lemma~4.4.4,
it is proven that
\[\chi_N(w,q)=\sum_{i=0}^k g_i(q)\bigl(q\,{d \over dq}\bigr)^i \chi_N(v,q) \]
with  $g_i(q)\in \C[[q]]$ for $i$, $\ldots$, $k$.
This implies that if $\chi_N(v,q)\in q^{h+1-c/24}\C[[q]]$ also 
$\chi_N(w,q)\in q^{h+1-c/24}\C[[q]]$. Since every element of $V$ is a sum of
elements as $w$, the theorem follows.
We also use that $\bigoplus_{k\leq t} V_{[k]}=\bigoplus_{k\leq t} V_{k}$
and $\bigoplus_{n\leq h} N_{[n]}=\bigoplus_{n\leq h} N_{n}$.
\phantom{xx} \eop

Let $G$ be a group of automorphisms acting on a \voa $V$. We say that a module $N$ of $V$ is
$G$-invariant, if there exists a central extension $ \widehat G$ of $G$  acting
on $N$ such that $g^{-1}Y_N(\bar g v,z)g=Y_N(v,z)$ for all $g\in   \widehat G$  and 
$v\in V$. Here, $ \bar g$ denotes the image of $g$ in $G$.

The following result is a generalization of~\cite{Ma-design}, Lemma~2.8.
\begin{theorem}\label{groupdesigns} 
Let $V$ be a \voa and $G$ be a compact Lie group of 
automorphisms of $V$. Let $X$ be a homogeneous subspace of a $G$-invariant module of $V$.
If  the minimal conformal weight of $V^G$ is larger or equal to $t+1$, then
$X$ is  a conformal $t$-design.
\end{theorem}

The analogous result that $t$-homogeneous permutation groups $G$ lead to 
block $t$-designs is trivial. An analogous result relating the invariants  
of a real representation of a finite group $G$ with spherical designs is 
due to Sobolev~\cite{So-quadratur}.

\smallskip

\pf For all $g\in  \widehat G$ and $v\in V$ we have $\tr|_X\,\o(\bar g v) =\tr|_X\,\o( v)$
as $\o(\bar g v)=g\o(v)g^{-1}$. Let $\tilde v$ be the average of the $\bar g v$, i.e.,  
$\tilde v=\frac{1}{|G|}\sum_{\bar g\in G}\bar gv$ for $G$ finite and $\tilde v=\int_G \bar g v\, \mu/\int_G 1\,\mu$
with $\mu$ the Haar measure on $ G$ in general. Then $\tr|_X\,\o( \tilde v)=\tr|_X\,\o( v)$.
Since $\tilde v\in V^G$ and the minimal weight of $V^G$ is larger  or equal to $t+1$
and the action of $G$ on $V_\omega$ is trivial, one
gets $\tilde v=\pi(v)$ and hence $\tr|_X\,\o(\pi(v))=\tr|_X\,\o( v)$ for 
$v\in V_s$, $0\leq s\leq t$. \eop

A $t$-design $X$ based on $V$ as in the previous theorem is trivial
as a conformal $t$-design based on $V^G$.

In~\cite{Ma-design}, a \voa $V$ 
is said to be {\it of class ${\cal S}^n$\/} if the minimal
weight of $V^G$ is larger or equal to $n$, where $G={\Aut}\,V$ is the automorphism group
of~$V$. (Matsuo assumes in addition $V_1=0$, but the definition
clearly works without this assumption.)

\begin{example}\label{a1}\rm
The lattice vertex operator algebra $V_{A_1}$ associated to the root lattice $A_1$ has
the complex Lie group ${\rm PSU}_2(\C)$ as automorphism group and $V_{A_1}$ is a
highest weight representation of level $1$ for the affine Kac-Moody algebra of type $A_1$.
The graded character of $V_{A_1}$ as a ${\rm PSU}_2(\C)$-module is
$$\chi_{V}=\Bigl(1+\sum_{n=1}^{\infty}(\lambda^n+\lambda^{-n})\,q^{n^2}\Bigr)\Bigm/
\Bigl(q^{1/24}\prod_{n=1}^{\infty}(1-q^n)\Bigr).$$
The characters of the irreducible ${\rm PSU}_2(\C)$-modules are $\sum_{k=-i}^i\lambda^k$ with 
non-negative integers $i$ and it follows
immediately that the graded multiplicity of the trivial ${\rm PSU}_2(\C)$-representation
in $V_{A_1}$ is $q^{-1/24}\left(\prod_{n=2}^{\infty}(1-q^n)\right)^{-1}$ which equals
the graded character of the Virasoro \voa of central charge~$1$. By restricting
${\rm PSU}_2(\C)$ to the compact group ${\rm PSU}_2(\R)$, Theorem~\ref{groupdesigns}
shows that  $V_{A_1}$ is of class ${\cal S}^t$ for all $t$ and 
the homogeneous subspaces $(V_{A_1})_n$ are conformal $t$-designs for all $t$.
By using the two-fold cover ${\rm SU}_2(\R)$ of ${\rm PSU}_2(\R)$ it follows that
the homogeneous subspaces of the irreducible $V_{A_1}$-module $V_{A_1+\frac{1}{\sqrt{2}}}$
of conformal weight $1/4$ are also conformal $t$-designs for all $t$.
\end{example}

\begin{example}\label{e8}\rm
The lattice vertex operator algebra $V_{E_8}$ associated to the root lattice $E_8$ has
the complex Lie group $E_8(\C)$ as automorphism group and $V_{E_8}$ is the unique 
highest weight representation of level $1$ for the affine Kac-Moody algebra of type $E_8$.
It can be deduced from the Weyl-Kac character formula that the minimal weight of
the fixed point \voa $V_{E_8}^{E_8(\C)}$ is~$8$. Hence, $V_{E_8}$ is of class ${\cal S}^7$ and 
the homogeneous subspaces $(V_{E_8})_n$ are conformal $7$-designs since $E_8(\C)$ can
be restricted to the compact group $E_8(\R)$ without changing the fixpoint \voa.
\end{example}

\begin{example}\label{moonshine}\rm
It was noted in~\cite{Ho-dr}, Sect.~5.1, that Borcherds' proof of the Moonshine conjectures 
implies that the minimal weight of $(V^{\natural})^{M}$,
where ${V^{\natural}}$ is the Moonshine module with the Monster $M$ as automorphism
group, is equal to $12$. Hence, $V^{\natural}$ is of class ${\cal S}^{11}$ and 
the homogeneous subspaces $V^{\natural}_n$ are conformal $11$-designs.
\end{example}

\begin{example}\label{shortmoonshine}\rm 
Using the proof of the generalized Moonshine conjecture in~\cite{Ho-babymoon}, it was
shown loc.~cit.~that the even part $V\!B^{\natural}_{(0)}$ of the odd Moonshine module 
$V\!B^{\natural}$ with the Baby Monster $B$ as automorphism group~\cite{Ho-shortaut} is of class 
${\cal S}^7$.
Hence, the homogeneous subspaces of $V\!B^{\natural}_{(0)}$ 
and of its two further $B$-invariant irreducible modules 
$V\!B^{\natural}_{(1)}$ and $V\!B^{\natural}_{(2)}$
are conformal $7$-designs.
\end{example}

We will give a more direct proof for Examples~\ref{moonshine} and~\ref{e8}
in the next section and for Example~\ref{shortmoonshine} after the next theorem.

\medskip 

Let  $V$ be a \voa of central charge $c$ and $N$ be a $V$-module  graded by $\Z+h$.
Let $U=V_{\omega'}$ be a Virasoro \subvoa of $V$ of central charge~$c'$ with  
Virasoro element $\omega'$, let $W$ be the commutant of $U$ in $V$ and assume  
that  $\omega^*=\omega-\omega'$ is a Virasoro element of $W$ of central charge $c^*=c-c'$. 
Assume also that as $U\otimes W$-module one has a finite decomposition 
$N=\bigoplus_{h'+h^*=h}\widetilde{M}(h')\otimes K(h^*)$, where $\widetilde{M}(h')$ is
a module of $V_{\omega'}$ generated by highest weight vectors of degree $h'$ and
the $K(h^*)$ are $W$-modules graded by $\Z+h^*$. For the homogeneous degree-$h$-part $X$ 
of $N$  one gets the  finite decomposition 
$$X=\bigoplus_{h'+h^*=h} \bigoplus_{\,k\in \Z_{\geq 0} } X'_{h',h'+k}\otimes X^*_{h^*,h^*-k}$$
with homogeneous parts
$ X'_{h',h'+k}=\widetilde{M}(h')_{h'+k}$  and  $X^*_{h^*,h^*-k}=K(h^*)_{h^*-k}$.

\smallskip

\noi{\bf Definition:} The homogeneous subspaces  $X^*_{h^*,h^*-k}$ of the 
$W$-modules $K(h^*)$ are called  the {\it derived parts of $X$ with respect 
to $\omega'$.\/} We also define
$S=\{h'\mid \widetilde{M}(h')\not =0\}$ and~$s^*=|S|$.

\smallskip

\begin{theorem}\label{derivedparts} Let $V$ be a \voa of central charge $c$ and
let $N$ be a $V$-module  graded by $\Z+h$. Assume that $\omega'$, $\omega-\omega'\in V_2$ 
generate two commuting Virasoro \subvoas of central charge~$c$ and $c^*=c-c'$, respectively.
Assume further that the homogeneous parts $N_n$ for $n\leq h$ are conformal designs
of type $t$-$(c,n)$ based on $V$. Then the derived parts $X^*_{h^*,h^*-k}$ of $X=N_h$ 
with respect to $\omega'$ are conformal designs of type $t^*$-$(c^*,h^*-k)$ based 
on the commutant $W$ of $V_{\omega'}$ in $V$ with $t^*=t+2-2s^*$.
\end{theorem}

\pf Let $\omega^*=\omega-\omega'$ be the Virasoro element of $W$. From the commutative diagram
$$\begin{array}{lclcl}
V_{\omega'} & \!\otimes \! & W & \hookrightarrow & V \\
\, \downarrow  \hbox{\small id} & & \, \downarrow \hbox{\small $\pi^*$} & &\, \downarrow \hbox{\small $\pi$} \\
V_{\omega'} &  \!\otimes \! & V_{\omega^*} & \rightarrow & V_{\omega}
\end{array}$$
where the upper horizontal arrow is the natural injection and the other arrows are the
projections onto the relevant Virasoro \subvoas it follows that if $w \in \ker \pi^*$,
then $u\otimes w \in \ker \pi$ for all $u\in V_{\omega'}$.

Let now $w\in \ker \pi^*\cap W_s$ where  $s\leq t^*$. 
Denote with $L'_{n}$ and $L^*_n$ the generators of the Virasoro algebra associated to~$\omega'$ 
and~$\omega^*$, respectively.
For $0\leq i \leq s^*-1$ we define
$v^i\in V_{s+2i}$ by $v^i=(L'_{-2})^i{\bf 1}\otimes w$.
For the graded traces one gets
$$\tr|_N\, \o(v^i) q^{L_0}=\sum_{h'\in S}(\tr|_{\widetilde{M}(h')}\o((L'_{-2})^i{\bf 1})\, q^{L_0'})\cdot 
                         (\tr|_{K(h-h')}\o(w)\, q^{L_0^*}).$$
The matrix valued power series 
$$(\tr|_{\widetilde{M}(h')}\o((L'_{-2})^i{\bf 1})\,q^{L_0'-h'}), \qquad h'\in S,\quad i=0,\, 1,\, \ldots,\,s^*-1,$$
has as constant term the matrix $(n_{h'}\,P_i(h'))$, where $n_{h'}$ is the dimension of the
lowest degree subspace $\widetilde{M}(h')_{h'}$ and $P_i(h')$ is a monic polynomial of degree~$i$ in $h'$.
This matrix is invertible since its determinant is of Vandermonde type 
and therefore the matrix valued power series is invertible, too.
Since the $N_n$ are conformal $t$-designs for $n\leq h$, the coefficients of 
$\tr|_N \o(v^i)\,q^{L_0}=\sum_{\Z+h}a_n^i\,q^n$ vanish for $n\leq h$ and for all $i$ as $s+2i\leq t^*+2(s^*-1)= t$.
Hence the coefficients of $\tr|_{K(h-h')}\o(w)\, q^{L_0^*+h'}=\sum_{k\in\Z+h-h'} b^h_k\,q^{h'+k}$ also vanish
for $k\leq h-h'=h^*$ and all $h'\in S$.  This implies that the $X^*_{h^*,h^*-k}$ are
conformal $t^*$-designs based on $W$. \phantom{xx}  \eop
 
\smallskip

We give two examples, both based on the Moonshine module $V^{\natural}$. As already
mentioned before, all its homogeneous subspaces $V^{\natural}_n$ are conformal $11$-designs.

\begin{example}\label{shortmoonshine2}\rm
By taking for $\omega'$ a Virasoro element of central charge $\frac{1}{2}$ one has the 
decomposition~(cf.~\cite{Ho-dr}, Sect.~4.1)
$$V^{\natural}=L_{1/2}(0)\otimes V\!B^{\natural}_{(0)}\, \oplus\,L_{1/2}(\frac{1}{2})\otimes V\!B^{\natural}_{(1)}\, 
\oplus\,L_{1/2}(\frac{1}{16})\otimes V\!B^{\natural}_{(2)},$$
where  the  $L_{1/2}(h')$ are the three irreducible Virasoro highest weight modules for
the Virasoro algebra of central charge $\frac{1}{2}$ and
the  $V\!B^{\natural}_{(i)}$ are the irreducible modules of $V\!B^{\natural}_{(0)}$, the even part
of the shorter Moonshine module~$V\!B^{\natural}$.
It follows from Theorem~\ref{derivedparts} that the homogeneous parts of the 
 $V\!B^{\natural}_{(i)}$ are conformal $7$-designs, i.e.,
the same result as already obtained from Theorem~\ref{groupdesigns}.
\end{example}

\begin{example}\label{laminated}\rm
The Leech lattice $\Lambda$ has a sublattice isomorphic to  $\Lambda_8 \oplus \Lambda_{16}$,
where $\Lambda_8$ is the $E_8$-root lattice rescaled by the factor $\sqrt{2}$ and 
$\Lambda_{16}$ is the Barnes-Wall lattice of rank~16. This sublattice
defines a \subvoa $V_{\Lambda_8}^+\otimes  V_{\Lambda_{16}}^+$ of 
$V^{\natural}$. The \voa $V_{\Lambda_8}^+$ was studied by Griess in~\cite{Gr-baby}
and has $O(10,2)^+$ as automorphism group. It can be easily seen 
from~\cite{DGH-virs} that  $V_{\Lambda_8}^+$ and $V_{\Lambda_{16}}^+$ are both
framed \voas. In fact, $V_{\Lambda_8}^+$ is isomorphic to the framed \voa $V_{\cal C}$ where the 
binary code ${\cal C}$ is equal to the Hammingcode ${\cal H}_{16}$ of length~$16$
and the code ${\cal D}$ is the zero code. The irreducible modules of a framed \voa with ${\cal D}=0$
are described in~\cite{Mi-rep}. One finds that $V_{{\cal H}_{16}}$ has $2^{10}$ modules 
$L_{\mu}$, $\mu \in \Xi$, having the  conformal weights  $0$, $\frac{1}{2}$ or $1$. 
Therefore one has a decomposition
\begin{equation}\label{e8decomp}
V^{\natural} = \bigoplus_{\mu \in \Xi} L_{\mu}\otimes K_{\mu} 
\end{equation}
of the Moonshine module into $V_{\Lambda_8}^+\otimes  V_{\Lambda_{16}}^+$-modules.
The graded traces $\chi_{L_{\mu}}(q)$ of the modules $L_{\mu}$ are easily computed. It turns
out that they  depend only on the conformal weight of $L_{\mu}$, i.e., they give
isomorphic $V_{\omega'}$-modules where $\omega' $ is the Virasoro element of $V_{\Lambda_8}^+$.
From decomposition~(\ref{e8decomp}),
we obtain now the following decomposition into $V_{\omega'}\otimes  V_{\Lambda_{16}}^+$-modules: 
\begin{equation}\label{omegadecomp}
V^{\natural} = \bigoplus_{m=0}^{\infty}L'(m) \otimes K(0)\, \oplus\,
                \bigoplus_{m=0}^{\infty}L''(m+\frac{1}{2}) \otimes K(\frac{3}{2})\, \oplus\,
                \bigoplus_{m=0}^{\infty}L'''(m+1) \otimes K(1),
\end{equation}
where the $L'(h)$, $L''(h)$, and $L'''(h)$ are direct sums of highest weight representations
of highest weight $h$ for the Virasoro algebra of central charge~$8$ associated to $\omega'$
and $K(0)= V_{\Lambda_{16}}^+$, $K(\frac{3}{2})$ and $K(1)$ are $V_{\Lambda_{16}}^+$-modules
of conformal weight $0$, $\frac{3}{2}$ and $1$, respectively.
Although Theorem~\ref{derivedparts} is not directly applicable, the argument given in its
proof shows that the homogeneous subspaces of $K(0)= V_{\Lambda_{16}}^+$, $K(\frac{3}{2})$, and $K(1)$
are conformal $7$-designs based on $V_{\Lambda_{16}}^+$. Moreover, one can switch the r\^ole of
$V_{\Lambda_8}^+$ and $  V_{\Lambda_{16}}^+$ in the preceding discussion and one obtains that
the  homogeneous subspaces of 
$V_{\Lambda_8}^+$ and the direct sum of the $V_{\Lambda_8}^+$-modules $L_{\mu}$ of 
conformal weight $\frac{1}{2}$ and $1$, respectively, are also conformal $7$-designs.
\end{example}

The application range of Theorem~\ref{derivedparts} is somewhat restricted 
since the minimal weight of the underlying \voa cannot be larger than $2$ and this
leaves only a limited set of examples of $t$-designs with $t\geq 6$ as
the classification Theorem~\ref{classfic7} of Section~\ref{classresult} shows. 

\medskip

One may ask if the \voas $V_L$ associated to even integral lattices~$L$ besides the 
root lattices $A_1$ and ${E_8}$ lead to  interesting conformal $t$-designs for larger values of~$t$. 
The next theorem shows that a necessary condition is that one starts with a  spherical $t$-design.
We recall that the irreducible modules of a lattice \voa associated to a lattice $L$
are parametrized by the elements of the discriminant group~$L^*/L$.

\begin{theorem}\label{latticevoadesign}
Let $V_{L}$ be the lattice \voa associated to an even  integral lattice 
$L\subset {\bf R}^k$ of rank $k$ and let $N=\sum_{i} V_{L+\lambda_i}$,
$\lambda_i\in L^*/L$, be a module of $V_L$.
If the degree-$n$-subspace $X=N_n$ of $N$ is a conformal $t$-design based on $V_L$
then the set of vectors of norm $2n$ in $\bigcup_i L+\lambda_i$ must form a 
spherical $t$-design. 
\end{theorem}
The proof uses the following result (\cite{DMN-latticetrace}, Th.~3):
\begin{prop}[Dong-Mason-Nagatomo]\label{DMN-harmonic}
Let $P$ be a homogeneous spherical harmonic polynomial on ${\bf R}^k$ and
let $V_L$ be the lattice \voa associated to an even integral lattice 
$L$ of rank~$k$. Then there exists a Virasoro highest weight vector
$v_P$ with the property that
\[
\chi_{V_L}(v_P,q)=
\Bigl(\sum_{x \in L}P(x)\,q^{(x,x)/2}\Bigr)\Bigm/
\Bigl( q^{1/24}\prod_{i=1}^{\infty}(1-q^i)\Bigr)^k.
\]
\end{prop}
The vector $v_P$ is given by
$v_P=P(h^1_{(-1)},\, \ldots,\, h^k_{(-1)}){\bf 1}$, where $\{h^1,\, \ldots,\,h^k\}$
is an orthonormal basis of ${\bf R}\otimes_{\bf Z} L\subset {(V_L)}_1$.
It can immediately be seen from its proof that the proposition remains valid
if one replaces $L$ by a coset $L+\lambda$ and $V_L$ 
by the $V_L$-module  $V_{L+\lambda}$.

\smallskip
{\bf Proof of Theorem~\ref{latticevoadesign}.} As explained in the introduction,
the set $\widetilde{L}_{2n}=\{x\in \bigcup_i L+\lambda_i \mid (x,x)=2n\}$ on a
 sphere in  ${\bf R}^k$ around $0$ 
is a spherical $k$-design if and only if $\sum_{x\in \widetilde{L}_{2n}}P(x)=0$ for all
harmonic polynomials $P$ homogeneous of degree $s$ with $1\leq s\leq t$. The result follows
now directly from the mentioned generalization of Proposition~\ref{DMN-harmonic}. \eop

It will follow from  Theorem~\ref{classfic7} part (b) of Section~4 that
for a conformal $t$-design based on a lattice \voa $V_L$ the largest $t$ 
one can hope for is $t=5$ if $L\not= A_1$, $E_8$.

\medskip

One of the main results of~\cite{Ma-design} are certain formulae
for traces of the form 
$$\tr|_{V_n}\o(v_1)\o(v_2)\cdots\o(v_k)$$
with $v_1$, $\ldots$, $v_k\in V_2$ and $V$ a \voa of class ${\cal S}^{2k}$,
$k\leq 5$, with $V_1=0$ and ${\rm Aut}\, V$ finite 
(Theorem 2.1 for $n=2$ and Theorem~5.1 for $k\leq 2$ and general $n$). 
For the Moonshine module $V^{\natural}$ and $n=2$ they were first obtained by 
S.~Norton~\cite{No-newformulae}. The proof given in~\cite{Ma-design}
remains valid, if one replaces the assumption that $V$ is
a \voa of class ${\cal S}^{2k}$ and ${\rm Aut}\, V$ is finite with the
assumption that $V_l$ for $l\leq n$ is a conformal $2k$-design based on $V$. 

One can study similar trace identities for conformal $t$-designs supported 
by a module of $V$, without the assumption $V_1=0$, and with 
$v_1$, $\ldots$, $v_k$ homogeneous elements of $V$ not necessarily in $V_2$.

\medskip

We end this section with an open problem:  It is known that there exist non-trivial block 
$t$-designs~\cite{Te-alltbd} and  spherical $t$-designs~\cite{SeZa-alltsd} 
for arbitrary large $t$ and arbitrary large length respectively dimension. 
The same result for block designs supported by self-orthogonal codes and 
for spherical designs supported by integral lattices
seems to be open (and less likely). The example of the lattice \voa $V_{A_1}$
(Example~\ref{a1}) shows that there exist non-trivial $t$-designs for arbitrary
large values of $t$. However, this case may be considered exceptional because the
central charge $c=1$ of $V_{A_1}$ is small. We ask therefore if there exist non-trivial
conformal $t$-designs for arbitrary large values of $t$ and arbitrary large 
central charge.


\section{Conformal designs associated to extremal \voas{}}\label{extremal}

A rational vertex operator algebra $V$ is called self-dual (other authors use
also the names holomorphic or meromorphic) if the only irreducible module
of $V$ up to isomorphism is $V$ itself. The central charge $c$ of a self-dual \voa is
of the form $c=8l$ where $l$ is a positive integer.
Its character $\chi_V$ is a weighted homogeneous polynomial of weight $c$ in
the polynomial ring over the rationals generated by 
the character of the self-dual lattice \voa $V_{E_8}$ associated to the $E_8$ lattice
(given the weight~$8$) and the character of the self-dual Moonshine module $V^{\natural}$
(given the weight~$24$); see~\cite{Ho-dr}, Thm.~2.1.2. Since
$$\chi_{V_{E_8}}= \sqrt[3]{j} = q^{-1/3}(1+248\,q + 4124\, q^2+ 34752\, q^3+\cdots )$$
and 
$$\chi_{V^{\natural}}= j-744=q^{-1}(1+196884\, q^2+ 21493760\,q^3+\cdots )$$
one can use alternatively $\sqrt[3]{j}$ (weight~$8$) and the constant function
$1$ (weight~$24$) as generators; see~\cite{Ho-dr}, Sect.~2.1.

In~\cite{Ho-dr}, extremal \voas were defined. A self-dual \voa $V$
of central charge $c$ is called {\it extremal\/} if its minimal
weight $\mu(V)$ satisfies $\mu(V)> \left[\frac{c}{24}\right]$, i.e., a
Virasoro primary highest weight vector of $V$ different from a multiple of the
vacuum has at least the conformal weight $\left[\frac{c}{24}\right]+1$.
It follows from the above description of the character of a self-dual \voa,
that the character of an extremal \voa $V$ has the form
$$\chi_V=\chi_{V_{\omega}}\cdot ( 1+ A_{k+1}\,q^{k+1}+  A_{k+2}\,q^{k+2}+\cdots),$$
with $k= \left[\frac{c}{24}\right]$ and constants $A_l$ independent of $V$.
It can be shown that
$A_{k+1}>0$ (see~\cite{Ho-dr}, Thm.~5.2.2), i.e., the minimal weight
of an extremal \voa is in fact equal to $\left[\frac{c}{24}\right]+1$.

\smallskip

\begin{theorem} \label{extremaldesign}
Let $V$ be an extremal \voa of central charge $c$. Then the
degree $n$ subspace $V_n$ of $V$ is a conformal $t$-design with 
$t=11$ for $c\equiv 0\!\!\pmod{24}$, $t=7$ for $c\equiv 8\!\!\pmod{24}$, 
and $t=3$ for $c\equiv 16\!\!\pmod{24}$.
\end{theorem}

\pf Let $v\in V_s$ be a Virasoro highest weight vector
of conformal weight $s$, where $0<s\leq t$. 
It follows from Zhu, 
\cite{Zhu-dr} Thm.~5.3.3, that 
$$\chi_{V}(v,q)=q^{-c/24} \sum_{n=0}^{\infty} \tr|_{V_n} \o(v)\,q^n$$
is a meromorphic modular form of weight $s$ for ${\rm PSL}_2(\Z)=
\langle S, T \rangle$ with character $\rho$ given by 
$\rho(S)=1$ and $\rho(T)=e^{-2\pi i c/24}$. 
Here, $S=\pm\bigl({\phantom{{-}}0\ 1 \atop {-}1 \ 0} \bigr)$ and
$T=\pm\left({1 \ 1 \atop 0\ 1} \right)$.
Since $V$ is assumed to be an extremal \voa, one has
\begin{equation}\label{tracezero}
\tr|_{V_n} \o(v)=\tr|_{(V_{\omega})_n} \o(v)=0
\end{equation}
for $n=0$,~$\ldots$,~$k$, where $k= \left[\frac{c}{24}\right]$. 
For the last equal sign in (\ref{tracezero}), one uses the skew-symmetry identity
$Y(v,z)u=e^{zL_{-1}}Y(u,-z)v$. This shows that for $u\in (V_{\omega})_n$ 
the product $\o(v)u\in V_n$ is contained in the Virasoro highest-weight module
generated by $v\in V_s$. Since $V$ is extremal, one has $s>k\geq n$ and therefore
$\o(v)u=0$.

Let $\Delta=q\prod_{n=1}^{\infty}(1-q^n)^{24}$ be the unique normalized cusp form 
of weight~$12$ for ${\rm PSL}_2({\Z})$. One has $\sqrt[24]{\Delta}=\eta$, the Dedekind
eta-function.
Since $\eta(-1/\tau)=\eta(\tau)$ and $\eta(\tau+1)=e^{2\pi i/24}\eta(\tau)$ it
follows that $\chi_V(v,q) \cdot \eta^{c}$ is a
holomorphic modular form for ${\rm PSL}_2(\Z)$ of weight $c/2+s$ and  trivial 
character for which the first $k+1$ coefficients of its $q$-expansion vanish.
Such a modular form is of the form $\Delta^{k+1}f$ for some
holomorphic modular form $f$ of weight 
$$c/2+s-12(k+1)= \cases{s-12, & for $c=24k$, \cr
s-8, & for $c=24k+8$, \cr s-4, & for $c=24k+16$. \cr}$$
Using the fact that there is no non-zero holomorphic modular form of 
negative weight, one concludes that $f=0$ if $s\leq t$ and $t$ as in
the theorem. This gives $\chi_V(v,q)=0$ and so $\tr|_{V_n} \o(v)=0$
for any $n$. The result follows now from Theorem~\ref{highestweight}\eop

\noi{\bf Remark:} Under the same assumption as in the theorem, our proof
shows that also $\tr|_{V_n} \o(v)=0$ for any Virasoro highest weight vector $v$
of conformal weight $s$ where
$s=13$, $14$, or $15$  for $c\equiv 0 \!\! \pmod{24}$, $s=9$, $10$ or $11$
for $c\equiv 8\!\! \pmod{24}$ and $s=5$, $6$, $7$ for
$c\equiv 16\!\! \pmod{24}$ if we use the fact that there is no 
non-zero holomorphic modular form of weight~$1$, $2$, or~$3$.

\smallskip

\begin{example}\label{extexa}\rm
Extremal \voas are known to exist for $c=8$, $16$, $24$, $32$ and $40$; 
see~\cite{Ho-dr}, Sect.~5.2. 

For the $c=8$, the only example is $V_{E_8}$
and we know already from Example~\ref{e8} that its homogeneous subspaces are
conformal $7$-designs. For $c=16$, the two self-dual \voas
$V_{E_8^2}$ and $V_{D_{16}^+}$ are both extremal and their homogeneous subspaces are
therefore conformal $3$-designs.

The known examples for $c=24$, $32$, $40$ are
$\Z_2$-orbifolds of lattice VOAs, where the lattice is an
even unimodular lattice of rank~$c$ without vectors of squared length~$2$,
i.e., an extremal lattice. For $c=24$, this gives only the Moonshine 
module $V^{\natural}$.
Using Theorem~\ref{groupdesigns}, we have already seen that its subspaces of 
fixed degree are conformal $11$-designs.  
For $c=32$, there exist at least $10^7$~\cite{king-uni32}
extremal even lattices. Our theorem shows in particular that the degree
subspace $V_2$ of the \hbox{$\Z_2$-orbifold} of the associated lattice \voa 
is a conformal $7$-design. In the case of the Barness-Wall lattice of rank~$32$,
it was observed by R.~Griess and the author that the automorphism group of the
\hbox{$\Z_2$-orbifold} \voa is likely $2^{27}{:}E_6(2)$ and one may use in this
case also Theorem~\ref{groupdesigns} to derive the conformal $7$-design property.
There are at least $12579$ extremal 
doubly-even codes~\cite{King-sd40} of length~$40$. Using orbifold constructions 
(cf.~\cite{DGH-virs}, Sect.~4) one sees that there are at least so many extremal 
even lattices of rank~$40$ and at least so many extremal \voas of central charge~$40$.
The homogeneous subspaces of those \voas, in particular $V_2$, 
are conformal $3$-designs. Characters of extremal \voas for small $c$ are
given in~\cite{Ho-dr}, Table~5.1. For $c=32$, one has $\dim\,V_2=139504$;
for $c=40$, one has $\dim\,V_2=20620$.

Since it is unknown if any extremal \voa of central charge~$48$, $72$,~$\ldots$
exists, we do not get currently any other conformal $11$-designs from our theorem
besides the ones from the Moonshine module. 
\end{example}

\smallskip

The minimal weight of an extremal \voa grows linearly with $c$ and therefore 
the conformal designs as in Theorem~\ref{extremaldesign} become trivial 
for $c\geq 264$, $176$, or $88$, if $c\equiv 0$,  $8$, 
or $16\!\!\!\pmod{24}$, respectively.

\medskip

One can ask for a  similar theorem  for \svoas. In~\cite{Ho-dr}, Sect.~5.3,
it has been shown that the minimal weight of a self-dual  \svoa
(under certain natural conditions called ``very nice'')  satisfies
$\mu(V)\leq \frac{1}{2}\bigl[\frac{c}{8}\bigr]+\frac{1}{2}$. Self-dual  \svoas
meeting that bound are called extremal.
Since that time, the analogous bounds for the minimal weight respectively~length
of even self-dual codes respectively unimodular lattices have been improved to the same bounds
as one knows for doubly-even self-dual codes and even unimodular lattices
(with the exception of codes of length $n\equiv 22 \!\!\pmod{24}$ and the shorter Leech
Lattice of rank~$23$).
So one may expect that for self-dual \svoas the analogous bound $\mu(V)\leq \bigl[\frac{c}{24}\bigr]+1+r$
holds with $r=\cases{\frac{1}{2}, & for $c \equiv 23\frac{1}{2} \!\pmod{24}$, \cr
0, & else.}$ For the case of the exceptional lengths $n\equiv 22 \!\!\pmod{24}$, it was
proven in~\cite{LaLa-codesign}, that codes meeting the improved bound lead to
block $3$-designs. The analogous result for \svoas would be that for 
a self-dual  \svoa $V$ of central charge $c=24k+23\frac{1}{2}$ and minimal
weight $\geq k+3/2$ the homogeneous subspaces of $V_{(0)}$-modules are conformal $7$-designs. 
For the proof, one has to analyze the singular part of the vector valued modular
functions associated to $V$.
However, in the case of lattices the analogous theorem only applies
to the shorter Leech lattice leading to spherical $7$-designs. 
The only known --- and likely the only --- example for \svoas would be the shorter 
Moonshine module $V\!B^{\natural}$ for which we have already shown in the previous
section that the homogeneous subspaces of $V\!B^{\natural}_{(0)}$-modules are 
conformal $7$-designs.


\section{Classification results}\label{classresult}

In this section, we investigate vertex operator algebras and super algebras
of minimal weight $m\leq 2$ whose degree-$m$-part form a conformal $6$- or $8$-design.

\subsection{Statement of results}

For conformal $6$-designs, we have the following classification result:
\begin{theorem}\label{classfic7}
Let $V$ be a simple \svoa of central charge $c$ and minimal weight $m$ 
and assume that $V$ has a real form such that the invariant bilinear form 
is positive-definite.
Denote with $V_{(0)}$ the even \subvoa of $V$.
If the degree-$m$-subspace $V_m$ is a conformal $6$-design, one has:
\begin{enumerate} 
\item[(a)] If  $m=\frac{1}{2}$, then  $V$ is isomorphic to the self-dual ``single fermion'' \svoa
$\VF\cong L_{1/2}(0)\oplus  L_{1/2}(\frac{1}{2})$ of central charge $1/2$.
\item[(b)] If $m=1$, then $V$ is isomorphic to the lattice  \voa
$V_{A_1}$ of central charge $1$ associated to the root lattice $A_1$ or the
lattice  \voa $V_{E_8}$ of central charge $8$ associated to the root lattice $E_8$.
\item[(c)] If $m=\frac{3}{2}$ and the additional assumption $\dim V_2>1$ holds,
the central charge of $V$ 
is either $c=16$ or $c=23\frac{1}{2}$. 
\item[(d)] For $m=2$ and the additional assumptions that $V$ is rational and
that for any $V_{(0)}$-module of conformal weight $h$ there exists a
$V_{(0)}$-module of the same conformal weight whose lowest degree subspace 
is a conformal $6$-design, it follows that $V$ is a \voa and
there are at most $11$ possible cases for the central charge.
The allowed values of $c$, $\dim V_2$ and conformal weights $h$ of a possible additional irreducible
$V_{(0)}$-modules are given in Table~\ref{listem2}. In particular, if $c\in\{24,\,32,\,40,\,
1496\}$, $V$ has to be self-dual.
\end{enumerate}
\end{theorem}
\begin{table}\caption{\small Possible central charges and conformal weights of
additional $V_{(0)}$-modules in case of Theorem~\ref{classfic7}~(d)}\label{listem2}
\small
$\begin{array}{|l|c|c|c|c|c|c|c|c|c|c|c|c|}
\hline
\hbox{cent.~charge $c$}\phantom{\frac{|^3}{|^3}}\!\!\!\!
 &   8 & 16 & \frac{808}{35} & \frac{47}{2} &
24 & 32 &  \frac{164}{5} & \frac{236}{7} &    \frac{242}{7} & 40 & 1496 
  \\ \hline
\dim V_2 &  156 & 2296 & 63428 & 96256 & 196884 &  139504  &  90118 &  63366 & 49291 & 20620 &  54836 \\ \hline
\hbox{conf.~weight $h$}\phantom{\frac{|^3}{|^3}}\!\!\!\!
  &   \frac{1}{2},1 &  1,\frac{3}{2} &
    \frac{103}{70},\frac{67}{35} &  \frac{3}{2},\frac{31}{16} & - & - &
    \frac{11}{5},\frac{12}{5} &  \frac{16}{7},\frac{17}{7} &
    \frac{67}{28},\frac{17}{7}  & - & - \\ \hline

\end{array}$
\end{table}

For all the cases in which examples are known, one gets in fact conformal $7$-designs. 
This theorem is analogous to similar results of Martinet~\cite{Ma-ladesign} on integral lattices
whose vectors of minimal norm $m\leq 4$ form a spherical $7$-design 
and of Lalaude-Labayle~\cite{LaLa-codesign} on binary self-orthogonal codes for which the
set of words of minimal weight $m\leq 8$ form a block $3$-design.

The additional assumption in Theorem~\ref{classfic7}~(d) that the lowest degree
subspaces of certain $V$-modules are also conformal $6$-designs seems quite strong.
It was introduced since it is necessary to apply Proposition~\ref{MODcond}. The analogue
of Proposition~\ref{MODcond} for codes and lattices holds without an analogous assumption,
so that this may be also the case for \voas. If one assumes only that for at least one
irreducible $V$-module with $h\not=0$ the lowest degree subspace is a conformal $6$-design, 
then one gets the same values for $c$ and $\dim V_2$ as in Table~\ref{listem2}, but no
further conditions on the allowed values for~$h$. If one drops the assumption on 
$V_{(0)}$-modules completely, all values of $c$ occurring in Lemma~\ref{diophant} are allowed.

By requiring that $V$ has a real form such that the restriction $(\,.\,,\,.\,)$ of 
the natural bilinear form to $V_2$ is positive-definite, it follows that the central charge $c$ of 
$V$ and the central charge $e$ of any \subvoa of $V$ with Virasoro element $\omega'\in V_2$ with
$L_1\omega'=0$ is positive since $0<(\omega',\omega')=\omega'_{(3)}\omega'+(L_1\omega')_{(2)}\omega'=
\omega'_{(3)}\omega'=e/2\cdot{\bf 1}$.

\smallskip

In addition, we obtain for conformal $8$-designs:
\begin{theorem}\label{classfic8}
\begin{enumerate} 
\item[(i)]
Let $V$ be a simple \svoa of central charge $c\not=\frac{1}{2}$, $1$ and minimal weight $m\leq 2$ 
with $\dim V_2>1$ having a real form such that the natural invariant bilinear form
is positive-definite.
If the weight-$m$-part $V_m$ is a conformal $8$-design,
then $V$ is a \svoa of central charge~$24$ and minimal weight~$2$ with
$\dim V_2=196884$. 
\item[(ii)]
If in addition, we assume 
that for any $V_{(0)}$-module of conformal weight $h$ there exists a
$V_{(0)}$-module of the same conformal weight whose lowest degree subspace 
is a conformal $8$-design,
then $V$ is a self-dual \voa  with the same
conformal character as the Moonshine module $V^{\natural}$.  
\end{enumerate} 
\end{theorem}

Analogous results for spherical $11$-designs and block $5$-designs characterizing the
Leech lattice and the Golay code can again be found in~\cite{Ma-ladesign} and~\cite{LaLa-codesign},
respectively.

\smallskip

Part (i) of Theorem~\ref{classfic8} for $m=2$ was proven by Matsuo~\cite{Ma-design}. 
He also considered the case (d) of Theorem~\ref{classfic7} with weaker assumptions resulting
in more candidates.

 
\medskip

In the next subsection, we will prove two Propositions which give
relations between $c$, $h$ and $\dim V_2$ for conformal $6$-designs.
In the subsequent four subsections, we will prove the four cases of Theorem~\ref{classfic7}. 
In subsection~\ref{8designs}, we will prove Theorem~\ref{classfic8}. 
In the final subsection~\ref{candidates}, we will discuss the  examples
respectively candidates of \svoas  occurring 
in Theorem~\ref{classfic7} and~\ref{classfic8} and compare our results
with the situation for codes and lattices.


\subsection{Conditions for $6$-designs}\label{cond6}

In this subsection, we assume that $V$ is a \voa of central charge $c$ and 
$\omega'$ and $\omega-\omega'$ are two elements in $V_2$ generating two commuting
Virasoro \voas $U_{\omega'}$ and $U_{\omega-\omega'}$ of central charge 
$e$ and $c-e$, respectively.
These two Virasoro \subvoas  are isomorphic to a quotient of the Verma module
quotient $M(x,0)/M(x,1)$ where $x=e$ or $c-e$.  We decompose now the subalgebra 
$U_{\omega'}\otimes U_{\omega-\omega'}$ of $V$
as a module for the Virasoro \subvoa $V_{\omega}$ of $V$.
\begin{lemma}\label{hwnumbers}
Let $\bigoplus_h\overline{M}(h)$ be the decomposition of $U_{\omega'}\otimes U_{\omega'-\omega'}$ into
isotypical components as a module for the Virasoro \voa $V_{\omega}$ of central charge $c$.
Assume that  $e$, $c-e\not=0$, $-\frac{22}{5}$, $-\frac{68}{7}$, $\frac{1}{2}$.  
Then the multiplicity $\mu_h$  of $M(h)$ in $\overline{ M}(h)$ is
$0$ for $h=1$, $3$, $5$, $7$ and is $1$, $1$, $1$, $2$ for $h=0$, $2$, $4$, $6$,
respectively. If also  $e$, $c-e\not=-\frac{46}{3}$, $-\frac{3}{5}$, then, in addition,
the multiplicity of $M(8)$ in $\overline{ M}(8)$ is $3$.

If $e=\frac{1}{2}$ and $U_{\omega'}$ is assumed to be simple and the other assumptions hold, then the
multiplicity of $M(6)$ in $\overline{ M}(6)$ is $1$, the multiplicity of 
$M(8)$ in $\overline{ M}(8)$ is $2$ and the other multiplicities are the same.
 \end{lemma}
\pf Let $U(x)$ be a Virasoro \voa of central charge $x$. For  $x\not=0$, 
$-\frac{22}{5}$, $-\frac{68}{7}$, $\frac{1}{2}$ (and  $-\frac{46}{3}$, $-\frac{3}{5}$)
one has  $\dim U(x)_n= \dim (M(x,0)/M(x,1))_n$
for $n\leq 7$ (respectively, $n\leq 8$) because the formula~(\ref{kacdet})
shows that under these conditions there are no singular vectors of degree $n$ in $M(x,0)/M(x,1)$. 
If we know that the Virasoro module $U_{\omega'}\otimes U_{\omega-\omega'}$ does not contain any additional 
singular vectors of degree $n\leq 6$ (respectively, $n\leq 8$) besides the one in $M(c,0)$, the 
result follows from
 $$q^c\chi_{U_{\omega'}\otimes U_{\omega-\omega'}}=
 \bigl(\prod_{n=2}^{\infty}(1-q^n)^{-1}\bigr)^2 +{\rm O}(q^9)\qquad\qquad\qquad\qquad\qquad\qquad\qquad\qquad$$
$$=\prod_{n=2}^{\infty}(1-q^n)^{-1} + (q^2+q^4+2\,q^6+3\,q^8) \prod_{n=1}^{\infty}(1-q^n)^{-1} + {\rm O}(q^9)
=\sum_{h} \mu_h  q^c \chi_{M(h)}. $$

In the case $e=\frac{1}{2}$ with simple $U_{\omega'}$ one has to form a quotient module 
of $M(1/2,0)/$ $M(1/2,1)$ by dividing out additional singular vectors and  one gets
$U_{\omega'}\cong  L(1/2,0)$ with character
$$\chi_{U_{\omega'}}=q^{-1/48}(1+q^2+q^3+2 q^4+2 q^5+3 q^6+3 q^7+5 q^8+5 q^9+\cdots).$$
Under the same assumption as before it follows
that the multiplicity of $M(h)$ in $\overline{M}(h)$ does not change for $h<6$ or $h=7$,  but
the multiplicity of $M(6)$ in $\overline{M}(6)$ is only~$1$ and the multiplicity of $M(8)$ 
in $\overline{M}(8)$ is~$2$.

The space of Virasoro highest weight vectors of degree $h$ in  $U_{\omega'}\otimes U_{\omega-\omega'}$ 
can also be computed explicitly for $n\leq 8$ by using the characterization $L_1v=L_2v=0$ 
for highest weight vectors $v\in \overline{ M}(h)$.
Its dimension equals indeed the given values for $\mu_h$.
\phantom{xxx}\eop

We choose highest weight vectors  $v^{(2)}\in \overline{ M}(2)$, 
$v^{(4)}\in \overline{ M}(4)$ and two linear independent highest weight vectors  
$v_a^{(6)}$, $v_b^{(6)}\in \overline{ M}(6)$.

We list here only $v^{(2)}$ and $v^{(4)}$ as the expressions for  $v_a^{(6)}$ and $v_b^{(6)}$ are
quite long:
\begin{eqnarray} \label{v2}
v^{(2)} & = & b_{-2}{\bf 1} - \frac{\left( c - e \right) }{e} \,a_{-2}{\bf 1},  \\
v^{(4)} & = &  \frac{-3}{5} \,b_{-4}{\bf 1} + b_{-2}^2{\bf 1} - 
   \frac{3\,\left( 22 + 5\,c - 5\,e \right) \,\left( c - e \right)}{5\,e\,\left( 22 + 5\,e \right)}\, a_{-4}{\bf 1}
 - \nonumber \\ \label{v4}
& &    \frac{2\,\left( 22 + 5\,c - 5\,e \right) }{5\,e} \,a_{-2}b_{-2}{\bf 1} + 
   \frac{\left( 22 + 5\,c - 5\,e \right) \,\left( c - e \right) \,
     }{e\,\left( 22 + 5\,e \right) } a_{-2}^2{\bf 1}.
\end{eqnarray}
Here, $a_{n}$ and $b_{n}$ denote the usual generators of the Virasoro algebras 
of  $U_{\omega'}$ and $U_{\omega-\omega'}$, respectively.
The expressions for $v_a^{(6)}$ and $v_b^{(6)}$ are not well-defined for $e=\frac{1}{2}$ or
$c-e=\frac{1}{2}$.

\smallskip

If $e=\frac{1}{2}$ and $c-e\not=\frac{1}{2}$,
the Virasoro module $M(1/2,0)/M(1/2,1)$ contains the singular vector
 $s_6=(-108\, a_{-6} - 264\, a_{-4}a_{-2} + 93\, a_{-3}^2 +  64 \,a_{-2}^3){\bf 1}$
of degree~$6$ and $(U_{\omega'})_6\cong (L_{1/2}(0))_6=(M(1/2,0)/M(1/2,1))_6/
\C\,s_6$.
A representative $v^{(6)}$ for a non-zero highest weight vector in $\overline{M}(6)$ is given by
\begin{eqnarray}
v^{(6)} \! &\!=\!&
  \bigl((-2/3)( 2373 + 1657\,c + 20\,c^2 ) \,b_{-6} - 
  (1/3) ( 12387 + 7301\,c + 112\,c^2 ) \,b_{-4}b_{-2} \nonumber \\
&& +   (1/3) ( 6051 + 854\,c + 70\,c^2 ) \,b_{-3}^2 
 +  ( 501 + 1099\,c ) \,b_{-2}^3 \nonumber \\
&& - (5/48)( -1 + 2\,c ) ( -5031 + 3195\,c + 1696\,c^2 + 140\,c^3 ) \,a_{-6} \nonumber \\
&& -  (11/3)( -5031 + 3195\,c + 1696\,c^2 + 140\,c^3 ) \, a_{-4}b_{-2} \nonumber \\
&&  -  (2/3) ( -5031 + 3195\,c + 1696\,c^2 + 140\,c^3 ) \,a_{-3}b_{-3}\nonumber \\
&& + (1/3)( -102555 + 89361\,c + 12970\,c^2 + 224\,c^3 ) \, a_{-2}b_{-4} \nonumber \\
&&  - 407( -129 + 115\,c + 14\,c^2 ) \, a_{-2}b_{-2}^2 \nonumber \\
&& - (7/24)( -1 + 2\,c )( -5031 + 3195\,c + 1696\,c^2 + 140\,c^3 ) \,a_{-4}a_{-2} \nonumber \\
&& + (35/192)( -1 + 2\,c )  ( -5031 + 3195\,c + 1696\,c^2 + 140\,c^3 ) \,a_{-3}^2 \nonumber \\
&& + 7( -5031 + 3195\,c + 1696\,c^2 + 140\,c^3 ) \,   a_{-2}^2b_{-2}\bigr){\bf 1}.
\end{eqnarray}

\smallskip

Assume that $W$ is a module of $V$ of conformal weight~$h$.
Using the associativity relation for a vertex operator algebra and its modules 
one can evaluate
the trace of $\o(v)$ for an element $v=a_{p_1}\cdots a_{p_k} b_{q_1}\cdots b_{q_l}{\bf 1}
 \in U_{\omega'}\otimes U_{\omega-\omega'}$ on the lowest degree part $W_h$.
One obtains
$\tr|_{W_h}\o(v)=\sum_{i=0}^n \alpha_i m_i^*$ for some $n$, where
$m_i^*=\tr|_{W_h} a_{0}^i$,  and  the $\alpha_i$ are explicit constants depending
on $c$, $e$ and~$h$. We also define $d^*=m_0^*=\dim W_h$. 
The traces of $v^{(2)}$, $v^{(4)}$, $v_a^{(6)}$,  $v_b^{(6)}$ and $v^{(6)}$
can now be computed.
Again, we list only the traces for $v^{(2)}$ and $v^{(4)}$ explicitly:
\begin{eqnarray}
\tr|_{W_h}\o(v^{(2)})&\! =\! & h\,d^* - \frac{c}{e}\, m_1^*, \\
\tr|_{W_h}\o(v^{(4)})&\! =\! & (\frac{h}{5}+h^2)\,d+  \frac{\left( 22 + 5\,c \right) \,
     \left( c - 2\,\left( e + 22\,h + 5\,e\,h \right)  \right) }{5\,e\,
     \left( 22 + 5\,e \right) }\, m_1^* \nonumber \\ 
&& + \frac{968 + 330\,c + 25\,c^2}{5\,e\, \left( 22 + 5\,e \right)}\, m_2^* . 
\end{eqnarray}

\begin{prop}\label{MODcond}
Let $V$ be a \voa of central charge $c\not= -24$, $-15$, $-\frac{44}{5}$, $\frac{34}{35}$, $1$
which contains elements $\omega'$, $\omega-\omega' \in V_2$ 
generating two commuting Virasoro \voas of central charge $e$ and $c-e$, respectively, 
with $e$, $c-e \not = -\frac{68}{7}$, $-\frac{22}{5}$, $0$. 
If $e=\frac{1}{2}$ or $e=c-\frac{1}{2}$, assume that the Virasoro \voa generated by $\omega'$
or $\omega-\omega'$, respectively, is simple and the other assumptions hold.
If there exists a module $W$ of $V$ of conformal weight~$h\not= 0$ 
such that the lowest degree part $W_h$ is a conformal $6$-design,
then
$$4 + 7\,c + c^2 - 124\,h - 31\,c\,h + 248\,h^2=0 .$$ 
\end{prop}

\pf
Assume first that $e\not =\frac{1}{2}$, $c-\frac{1}{2}$. 
In this case the conditions on the central charges
guarantee that the Virasoro \voas  generated by
 $\omega'$ and $\omega-\omega'$ are isomorphic to $M(x,0)/M(x,1)$
up to degree~$6$, where $x=e$ or $c-e$, respectively. 
By assumption, the degree~$h$ subspace
$W_h$ is a conformal $6$-design. Therefore,
the highest weight vectors $v^{(2)}$, $v^{(4)}$, $v_a^{(6)}$ and  $v_b^{(6)}$
give the trace identities
\begin{equation}\label{modtraces}
\tr|_{W_h}\o(v^{(2)})=\tr|_{W_h}\o(v^{(4)})=\tr|_{W_h}\o(v_a^{(6)})=\tr|_{W_h}\o(v_b^{(6)})=0.
\end{equation}
By evaluating the traces, 
the equations~(\ref{modtraces}) lead to  a homogeneous system of linear equations 
for $d^*=m_0^*$, $m_1^*$, $m_2^*$
and $m_3^*$. Since $d^*>0$, the system has to be singular and its determinant
$$\Delta=\frac{15\,(c+24)(c+15)(c+\frac{44}{5})(c-\frac{34}{35})( 4 + 7\,c + c^2 - 124\,h - 31\,c\,h + 248\,h^2 )(c-e)h}
{8\,(e+\frac{68}{7})(e+\frac{22}{5})^2\,e^2\,(e-\frac{1}{2})}  $$
has to vanish. The proposition follows in this first case.

\smallskip
If $e =\frac{1}{2}$ or $e=c-\frac{1}{2}$, we can, without loss of generality, assume that 
$e=\frac{1}{2}$ and  $c-e\not=\frac{1}{2}$ because otherwise one can switch the role of 
$\omega'$ and $\omega-\omega'$ since $c\not = 1$.

Let $d_l^*$ be the dimension of the $a_0=\omega'_{(1)}$ eigenspace 
for the eigenvalue $l$. 
One has $m_0^*=d^*=\dim W_h = d_0^* +d_{1/2}^*+d_{1/16}^*$ and
$m_i^*=d_{1/2}^*\cdot (1/2)^i+d_{1/16}^* \cdot (1/16)^i$, for $i=1$, $2$, $3$.
The homogeneous system 
\begin{equation}\label{modtraces2}
\tr|_{W_h} \o(v^{(2)})=\tr|_{W_h} \o(v^{(4)})=\tr|_{W_h} \o(v^{(6)})=0
\end{equation}
of linear equations for $d_0^*$, $d_{1/2}^*$, $d_{1/16}^*$
has to be singular since $d^*>0$. Hence the determinant
$$\Delta= \frac{15(  c + 24 ) ( c+ 15 ) 
     ( c+\frac{44}{5} ) ( c-\frac{34}{35} ) \,
     ( 4 + 7\,c + c^2 - 124\,h - 31\,c\,h + 248\,h^2 )\, h }
{512 }  $$
of the system must vanish. The proposition follows also in this case.
\eop

\medskip

Using the associativity relation of a vertex operator algebra one can also evaluate
the trace of $\o(v)$ for $v=a_{p_1}\cdots a_{p_k} b_{q_1}\cdots b_{q_l}{\bf 1}$ on $V_2$.
One obtains now
$\tr|_{V_2}\o(v)=\beta+\sum_{i=0}^n \alpha_i m_i$ for some $n$, with 
$m_i=\tr|_{V_2} a_{0}^i$,  and  the $\alpha_i$ and $\beta$ are explicit constants depending
on $c$ and $e$. We also define  $d=m_0=\dim V_2$. The  traces of $v^{(2)}$, $v^{(4)}$, 
$v_a^{(6)}$, $v_b^{(6)}$ and $v^{(6)}$ can now be computed.
Again, we list only the traces for  $v^{(2)}$ and  $v^{(4)}$ explicitly:
\begin{eqnarray}
\tr|_{V_2}\o(v^{(2)})&\! =\! & 2\,d - \frac{c}{e}\,m_1 , \\
\tr|_{V_2}\o(v^{(4)})&\! =\! &  \frac{\left( 44 + 5\,c \right) \,\left( c - e \right) }{22 + 5\,e}+
  \frac{22}{5}\,d + 
   \frac{\left( 22 + 5\,c \right) \,
      \left( c - 22\,\left( 4 + e \right)  \right) }{5\,e\,
      \left( 22 + 5\,e \right) }\, m_1  \nonumber  \\ \label{traceid}
 &\!\! & \qquad\qquad \qquad \qquad\qquad \qquad 
                 +  \frac{968 + 330\,c + 25\,c^2}{5\,e\,  \left( 22 + 5\,e \right) }\, m_2 . 
\end{eqnarray}

\begin{prop}[Matsuo]\label{m2VOAcond}
Let $V$ be a \voa of central charge
$c\not= -24$, $-15$, $-\frac{44}{5}$, $\frac{34}{35}$, $\frac{55\pm \sqrt{33}}{2}$, $1$
with $V_1=0$ which contains elements $\omega'$, $\omega-\omega' \in V_2$ 
generating two commuting Virasoro \voas of central charge $e$ and $c-e$, respectively, with 
$e$, $c-e \not = -\frac{68}{7}$, $-\frac{22}{5}$, $0$. 
If $e=\frac{1}{2}$ or $e=c-\frac{1}{2}$, assume that the Virasoro \voa generated by $\omega'$
or $\omega-\omega'$, respectively, is simple and the other assumptions hold.
If $V_2$ forms a conformal $6$-design then 
$$d=\frac{c\,(2388+955\,c+70\,c^2)}{2(748-55\,c+c^2)}.$$ 
\end{prop}
\pf
Assume first $e\not =\frac{1}{2}$, $c-\frac{1}{2}$.
As in the previous proposition, the conditions on the central charges
guarantee that $v^{(2)}$, $v^{(4)}$, $v_a^{(6)}$ and  $v_b^{(6)}$
are defined. By assumption, $V_2$ is a conformal $6$-design. Therefore one has for 
the highest weight vectors $v^{(2)}$, $v^{(4)}$, $v_a^{(6)}$ and  $v_b^{(6)}$
the trace identities
\begin{equation}\label{min2traces}
\tr|_{V_2}\o(v^{(2)})=\tr|_{V_2}\o(v^{(4)})=\tr|_{V_2}\o(v_a^{(6)})=\tr|_{V_2}\o(v_b^{(6)})=0.
\end{equation}
The equations~(\ref{min2traces}) form a system of linear equations for $d$, $m_1$, $m_2$
and $m_3$. This system is non-singular since its determinant
$$\Delta=\frac{15(c+24)(c+15)(c+\frac{44}{5})(c-\frac{34}{35})(c-\frac{55-\sqrt{33}}{2})(c-\frac{55+\sqrt{33}}{2})(c-e)}
{4\,(e+\frac{68}{7})(e+\frac{22}{5})^2\,e^2\,(e-\frac{1}{2})}  $$
is not $0$. In this case there is a unique solution for $d$, $m_1$, $m_2$
and $m_3$. The solution for $d$ is given in the proposition and does not depend on $e$.

\smallskip
If $e =\frac{1}{2}$ or $e=c-\frac{1}{2}$, we can, without loss of generality, assume that 
$e=\frac{1}{2}$ and  $c-e\not=\frac{1}{2}$ because otherwise one can switch the role of 
$\omega'$ and $\omega-\omega'$ since $c\not= 1$.

Similar as in the proof of Proposition~\ref{MODcond},
let $d_l$ be the dimension of the $a_0=\omega'_{(1)}$ eigenspace 
for the eigenvalue $l$. 
One has $m_0=d=\dim V_2 =1+d_0 +d_{1/2}+d_{1/16}$ and
$m_i=2^i+d_{1/2}\cdot (1/2)^i+d_{1/16} \cdot (1/16)^i$, for $i=1$, $2$, $3$.
The system 
$$\tr|_{V_2} \o(v^{(2)})=\tr|_{V_2} \o(v^{(4)})= \tr|_{V_2} \o(v^{(6)})=0$$ 
of linear equations for $d_0$, $d_{1/2}$, $d_{1/16}$
is non-singular since its determinant
$$\Delta=  \frac{15(c+ 24 ) (c+ 15 ) (c+ \frac{44}{5}) (c- \frac{34}{35})
    (c-\frac{55-\sqrt{33}}{2})(c-\frac{55+\sqrt{33}}{2})  }
                {256} $$
is not $0$.
The solution for $d$ is again the one given in the proposition.
\eop

The last proposition was first obtained by Matsuo~\cite{Ma-design} using a different argument.

\medskip

To apply Proposition~\ref{MODcond} and~\ref{m2VOAcond} to a \voa as in Theorem~\ref{classfic7}
and~\ref{classfic8}, one has to find suitable elements $\omega'\in V_2$ such that $\omega'$
and $\omega-\omega'$ generate commuting Virasoro \voas.
By Theorem 5.1 of~\cite{FreZhu}, this is the case if $\omega'$ generates a
Virasoro \voa and $L_1\omega'=0$. The last condition will be automatically
satisfied if either $V_1=0$ or $\omega'$ is the Virasoro element of an affine
Kac-Moody \voa $U<V$ or a Clifford \svoa $U<V$: In the case of a Kac-Moody
\subvoa, $\omega'$ is given by the Sugawara expression 
$\omega'=\sum_{i}\alpha_i\,a^i_{(-1)}a^i$, where the $a^i$ form an orthonormal
basis of a non-degenerated invariant bilinear form on $U_1$. For a Virasoro highest weight
vector $v$ of conformal weight $h$ one has $[L_m,v_{n}]=((h-1)m-n)v_{m+n}$ with
$v_k=v_{(k-h+1)}$ (cf.~\cite{Kac-VOA}, Cor.~4.10). For $u\in U_1$ 
one gets therefore 
$$L_1u_{(-1)}u=[L_1,u_{-1}]u+u_{-1}L_1u=u_0u=0$$
and the claim follows. Similarly, for $b\in U_{1/2}$ and $\omega'=\frac{1}{2}b_{-3/2}b_{-1/2}{\bf 1}=
\frac{1}{2}b_{-3/2}b$ one gets 
$$L_1b_{-3/2}b=[L_1,b_{-3/2}]b+b_{-3/2}L_1b=  b_{-1/2}b=0$$
and the claim follows also in the case of Clifford \svoas.


\subsection{Minimal weight $\frac{1}{2}$}\label{min0.5}  

Let $V$ be a \svoa as in part (a) of Theorem~\ref{classfic7} and let $V_{(0)}$ be
the even \subvoa.
The \subsvoa of~$V$ generated by the degree-$1/2$-part $V_{1/2}$ is known to be 
isomorphic to the \svoa $V_{\rm Fermi}^{\otimes d^*}$, where $d^*=\dim V_{1/2}$ and 
$V_{\rm Fermi}\cong L_{1/2}(0)\oplus L_{1/2}(\frac{1}{2})$
is the so-called single fermion \svoa  of central charge $1/2$.
By assumption one has $\mu(V)=\frac{1}{2}$ and hence $d^*=\dim V_{1/2}\geq 1$.

\smallskip

Recall that the central charge $c$ of $V$ is positive and assume 
first $c\not=\frac{1}{2}$, $\frac{34}{35}$,~$1$. The Virasoro element $\omega'$
of a subalgebra $V_{\rm Fermi}\subset V_{\rm Fermi}^{\otimes d^*} $
generates a Virasoro algebra of central charge $e=\frac{1}{2}$.
Since $V_{1/2}$ is assumed to be a conformal \hbox{$6$-design},
one can apply now Proposition~\ref{MODcond} to the $V_{(0)}$-module $V_{(1)}$ of 
conformal weight $h=\frac{1}{2}$.
For $h=\frac{1}{2}$, the proposition gives for the central charge
either $c=\frac{1}{2}$ or $c=8$.

If $c=8$, one gets  as in the proof of Proposition~\ref{MODcond} by  using the equations~(\ref{modtraces})
that
\begin{equation}\label{mstars}
 d_0= {\frac{255\, d^*}{496}},\ \ \ \ \ \ 
 d_{1/2}=  {\frac{d^*}{496}},\ \ \ \ \ \ 
 d_{1/16}=   {\frac{15\,d^*}{31}}. 
\end{equation}
This implies that $496|d^*$ and so $d^*\geq 496$, contradicting $d^*/2\leq c$
which follows from $V_{\rm Fermi}^{\otimes d^*}\subset V$ and the fact that the Virasoro
element $\omega^*$ of $V_{\rm Fermi}^{\otimes d^*}$ satisfies $L_1\omega^*=0$.

\medskip

If $c=1$, one finds again a Virasoro element $\omega'$ of central charge $\frac{1}{2}$
for a subalgebra  $V_{\rm Fermi}\subset V$ and $\omega-\omega'$ is the Virasoro
element for another central charge $\frac{1}{2}$ Virasoro algebra.
Thus $V$ contains a subalgebra isomorphic to $L_{1/2}(0)^{\otimes 2}$. 
The only irreducible modules of $L_{1/2}(0)^{\otimes 2}$ of conformal weight $1/2$
are  $L_{1/2}(0)\otimes L_{1/2}(\frac{1}{2})$ and  $L_{1/2}(\frac{1}{2})\otimes L_{1/2}(0)$.
Thus  $d_{1/16}=0$. The two vectors $v^{(2)}$ and $v^{(4)}$ are also well-defined
if $e=c-e=\frac{1}{2}$. Together with the equations
$$\tr|_{W_h}\o(v^{(2)})=\tr|_{W_h}\o(v^{(4)})=0$$
one gets $d_0=d_{1/2}=0$ and so $d^*=0$, contradicting $d^*>0$.

The central charge $c=\frac{34}{35}$ belongs to the unitary minimal series
\hbox{$c=1-6/(n(n+1))$}, $n=3$, $4$, $\ldots$, for $n=14$ and hence $L_{34/35}(0)\subset V$.
The only irreducible $L_{34/35}(0)$ modules of half-integral highest weight are
$L_{34/35}(0)$, $L_{34/35}(2)$, $L_{34/35}(23)$ and $L_{34/35}(39)$ of
conformal weight $0$, $2$, $23$~and $39$, respectively.
The simple current module $L_{34/35}(39)$ is 
the only module which can be used to extend $L_{34/35}(0)$ to a simple \svoa $V$. 
Hence the minimal weight of $V$ would be at least $39$, a contradiction.

For $c=\frac{1}{2}$ one has $L_{1/2}(0)\subset V$. 
The two further irreducible modules of $L_{1/2}(0)$ have the conformal 
weights~$\frac{1}{2}$ and~$\frac{1}{16}$. 
It follows $V\cong L_{1/2}(0)\oplus L_{1/2}(\frac{1}{2})\cong V_{\rm Fermi}$ since
$L_{1/2}(\frac{1}{2})$ is a simple current for $L_{1/2}(0)$ and thus any non-trivial 
simple \svoa extension of $L_{1/2}(0)$ contains the module $L_{1/2}(\frac{1}{2})$ 
with multiplicity~$1$ and is unique up to isomorphism.

\smallskip

This finishes the proof of Theorem~\ref{classfic7} part (a).


\subsection{Minimal weight $1$}\label{min1}

We begin with the following identity:
\begin{lemma}\label{v1traceid}
Let $V$ be a \voa and $a^1$, $a^2$, $\ldots$, $a^l$  be elements of $V_1$. Then
$$
\tr|_{V_1}\o(a^1_{(-1)}a^2_{(-1)}\cdots a^l_{(-1)}{\bf 1}) =
\tr|_{V_1}a^l_{(0)}a^{l-1}_{(0)}\cdots a^1_{(0)}  \qquad \qquad \qquad \qquad \qquad$$
$$ \qquad +\ \tr|_{V_1}a^1_{(-1)}a^l_{(1)}a^{l-1}_{(0)}\cdots a^3_{(0)} a^2_{(0)}+ \\
 \tr|_{V_1}a^2_{(-1)}a^l_{(1)}a^{l-1}_{(0)}\cdots a^3_{(0)}a^1_{(0)}
+ \cdots  $$
$$ \qquad \ \ +\ \tr|_{V_1}a^{l-1}_{(-1)}a^l_{(1)}a^{l-2}_{(0)}\cdots a^2_{(0)}a^1_{(0)} 
 \ + \ \tr|_{V_1}a^{l}_{(-1)}a^{l-1}_{(1)}a^{l-2}_{(0)}\cdots a^2_{(0)}a^1_{(0)}.\qquad $$
\end{lemma}
\pf The result follows from application of the associativity  relation by induction.
\phantom{xx}\eop 
For $l=4$ the identity can be found in~\cite{Hurley-dr} in the proof of Lemma~5.2. 

\smallskip

Let $V$ now be a \svoa as in Theorem~\ref{classfic7}, part (b).
We will use that under our assumptions $V_1$ is a reductive Lie algebra
under the product $x_{(0)}y$ for $x$, $y \in V_1$
and that the \subvoa $\langle V_1 \rangle$ generated by $V_1$ 
is an integrable highest weight representation of the associated affine 
Kac-Moody Lie algebra (cf.~\cite{Ho-dr} for \voas with positive-definite 
bilinear form on a real form; see~\cite{DoMa-reductive,DoMa-integrable}
for a result which applies to rational vertex operator algebras). 

\begin{lemma}\label{samecharge}
The central charge of the \subvoa $\langle V_1 \rangle$ spanned by $V_1$ 
equals the central charge of $V$.
\end{lemma}
\pf Let $\omega'$ be the Virasoro element of $\langle V_1 \rangle$ and
denote the central charge of $\langle V_1 \rangle$ by $e$.
One has $e>0$ as $V_1\not=0$. Assume that $e<c$. Then the vector
$v_2$ given in~(\ref{v2}) is a Virasoro highest weight vector of $V$
and equation~(\ref{modtraces}) together with the $6$-design property of $V_1$  gives
$\tr|_{V_1}\o(v_2)=\dim V_1 -\frac{c}{e} \dim V_1=0$ and so $e=c$, a contradiction.
\phantom{xxxxxxx}\eop

\smallskip

As reductive Lie algebra, $V_1$ can be decomposed into a Lie algebra direct sum  
\[V_1={\bf t} \oplus {\bf g}_1 \oplus \cdots \oplus {\bf g}_m \]
with abelian Lie algebra  ${\bf t}$ and simple Lie algebras 
${\bf g}_1$, $\ldots$, ${\bf g}_m$. Let $\h$ be a Cartan subalgebra of $V_1$ 
and let $\{h^1,\,\ldots,\,h^k\}$ be an orthonormal basis
of $\h$ with respect to the invariant form $\langle\,,\,\rangle $ on $V_1$ induced
by the canonical invariant bilinear form on $V$, i.e.,
$\langle x,y \rangle {\bf 1}=x_{(-1)}y$ for $x$, $y\in V_1$. The form
$\langle \,.\,,\,.\, \rangle$
is an  orthogonal sum of a non-degenerated form on~${\bf t}$ and some
nonzero multiples of the Killing forms on  each of the 
simple factors of $V_1$. 

\begin{lemma}\label{root}
Let $P(x_1,\,\ldots,\,x_k)$ be a complex harmonic polynomial
in variables $x_1$, $\ldots$, $x_k$. Then
$$\tr|_{V_1}\o(P(h^1_{(-1)},\,\ldots,\,h^k_{(-1)}){\bf 1})=
        \sum_{\alpha\in \Phi} P(\alpha(h^1),\,\ldots,\,\alpha(h^k)),$$
where $\Phi \subset \h^*$ is the root system of $V_1$ corresponding to $\h$.
\end{lemma}
\pf 
Let $\h \oplus \bigoplus_{\alpha\in\Phi} L_{\alpha}$ be the root space decomposition of 
$V_1$ with respect to $\h$. Each $L_{\alpha}$ is one-dimensional  and $[h,a]=\alpha(h)a$
for $a\in L_{\alpha}$ and $h\in \h$. Hence
\[\tr|_{V_1}  P(h^1_{(0)},\,\ldots,\,h^k_{(0)})=
        \sum_{\alpha\in \Phi} P(\alpha(h^1),\,\ldots,\,\alpha(h^k)).\]
We claim that the trace of a monomial $T=h^{i_1}_{(-1)}h^{i_2}_{(1)}h^{i_3}_{(0)}\cdots h^{i_l}_{(0)}$
on $V_1$ where $i_1$, $\ldots$, $i_l \in\{1,\,\ldots\,k\}$ is $0$ for $l\geq 3$: Indeed,
for an element $v\in L_{\alpha}$ one has 
$Tv=\lambda \langle h^{i_2},v \rangle h^{i_1} \in {\bf h}$ with a constant $\lambda$; 
for $v\in {\bf h}$ one has $Tv=0$ if $l\geq 3$. 

For a harmonic polynomial $P=\sum_{i,\,j=1}^k a_{ij} x_ix_j$ of degree~$2$ we have
\[ \tr|_{V_1} \sum_{i,\,j=1}^k a_{ij} \cdot (h^{i}_{(-1)}h^{j}_{(1)}) 
        = \tr|_{\h}  \sum_{i,\,j=1}^k a_{ij} \cdot(h^{i}_{(-1)}h^{j}_{(1)})\qquad\qquad\qquad \qquad\qquad \]
\[ \qquad\qquad\qquad\qquad\qquad\qquad\qquad
=\sum_{i,\,j,\, m=1}^k a_{ij} \langle h^m,h^i\rangle\langle h^j,h^m\rangle
 = \sum_{i=1}^k a_{ii} = 0, \]
where the last equality holds because $P$ is assumed to be harmonic. Similarly,
$$\tr|_{V_1} \sum_{i,\,j=1}^k a_{ij} \cdot (h^{j}_{(-1)}h^{i}_{(1)})=0.$$ 

The result follows now from Lemma~\ref{v1traceid}
since the cases $l=0$ and $l=1$ are clear.\eop 

For $P=x_1^4-6x_1^2x_2^2+x_2^4$, the previous Lemma  can be found in~\cite{Hurley-dr}, 
Sect.~6 and~7.

\begin{lemma}\label{cartan-hw}
If $P(x_1,\,\ldots,\,x_k)$ is a complex harmonic polynomial, then 
the  vector $v_P=P(h^1_{(-1)},\,\ldots,\,h^k_{(-1)}){\bf 1}$
is a highest weight vector for the Virasoro algebra of $V$.
\end{lemma}
\pf Let $V({\bf h})$ be the Heisenberg \voa of central charge $k$ generated by ${\bf h}$
and let $\omega'$ be its Virasoro element. For the Virasoro algebra associated to 
$\omega'$ it was proven in~\cite{DMN-latticetrace}, Lemma~5.1.12, that $v_P$ is a highest weight
vector. Since $\omega''=\omega-\omega'$ lies in the commutant of $V({\bf h})$, one
has $\omega_{(2)}''v_p=0$ and $\omega_{(3)}''v_p=0$ and therefore
also $L_1 v_P=\omega_{(2)}v_p=0$ and $L_2 v_P=\omega_{(3)}v_p=0$, i.e., $v_P$ is also a highest 
weight vector for the Virasoro algebra of $V$.
\eop

\begin{lemma}[Hurley~\cite{Hurley-dr} and~\cite{Hurley-trace}, Lemma~6.2]
For  the root system $\Phi$ of a nonsimple and nonabelian reductive Lie algebra ${\bf g}$
of rank~$n$ there exists a complex harmonic polynomial $Q$ of degree $4$ in $n$ variables such that
$$  \sum_{\alpha\in \Phi} Q(\alpha(h^1),\,\ldots,\,\alpha(h^n))\not =0$$
for an  orthogonal basis $\{h^1,\,\ldots,\,h^n\}$ of a Cartan algebra $\h$.
\end{lemma}
\pf The rank~$n$ of ${\bf g}$ is at least $2$.
Let $Q$ be the above mentioned harmonic polynomial $x_1^4-6x_1^2x_2^2+x_2^4$. 
The root system $\Phi$ is the union $\Phi_1\cup\ldots\cup\Phi_m$ of the root systems
of the either simple or abelian components $\k_i$ ($i=1$, $\ldots$, $m$) of ${\bf g}$.
(For the maybe existing abelian component, the root system is empty.) We choose an orthogonal
basis  $\{h^1,\,\ldots,\,h^n\}$ of the real subspace of $\h$ such that $h^1$ and $h^2$ lie 
in the Cartan subalgebras
of the components $\k_1$ and $\k_2$, respectively. Then $\alpha(h^1) =0$ for any 
$\alpha \not \in \Phi_1$ and $\alpha(h^2) =0$ for any $\alpha \not \in \Phi_2$. 
We obtain therefore
\begin{equation}\label{rootsum}
\sum_{\alpha\in \Phi} Q(\alpha(h^1),\,\ldots,\,\alpha(h^n))=\sum_{\alpha\in \Phi_1}\alpha(h^1)^4+
\sum_{\alpha\in \Phi_2}\alpha(h^2)^4.
\end{equation}
Both sums on the right hand side of equation (\ref{rootsum}) are real and nonnegative since
we have assumed that the $h^i$ are in the real subspace of $\h$. Furthermore, at least one 
of the two root systems $\Phi_1$ and $\Phi_2$ is nonempty and spans the dual Cartan algebra 
of the corresponding component. Thus we can find a root $\alpha \in  \Phi_1 \cup \Phi_2$ 
such that either $\alpha(h_1)\not=0$
or $\alpha(h_2)\not=0$ and so $\sum_{\alpha\in \Phi} Q(\alpha(h^1),\,\ldots,\,\alpha(h^n))\not =0$.
\eop

\begin{lemma}
For the  root system $\Phi$ of a simple Lie algebra of rank $n$ not of type $A_1$ or $E_8$ 
there exists a  complex harmonic polynomial $Q$ of degree less than or equal to $6$ in $n$ variables such that
$$  \sum_{\alpha\in \Phi} Q(\alpha(h^1),\,\ldots,\,\alpha(h^n))\not =0.$$
\end{lemma}
\pf
One easily checks that for $y\in {\bf R}^n$ the following two polynomials of degree~$4$ and~$6$ 
are harmonic:
\begin{eqnarray*}
R_4 &= & (x,y)^4-\frac{6}{4+n} (x,y)^2||x||^2||y||^2+\frac{3}{8+6n+n^2}||x||^4||y||^4,\\ 
R_6 &= & (x,y)^6-\frac{15}{8+n} (x,y)^4||x||^2||y||^2+\frac{45}{48+14n+n^2}(x,y)^2||x||^4||y||^4\\
& &  \qquad\qquad   {}-\frac{15}{192+104n+18n^2+n^3}||x||^6||y||^6. \\ 
\end{eqnarray*}
\vspace{-1mm}
Here, $x=(x_1,\,\ldots,\, x_n)$, the standard scalar 
product on ${\bf R}^n$ is denoted by $(\,.\,,\,.\,)$, and we let $||x||^2=(x,x)$. 
(These polynomials can be obtained from Gegenbauer polynomials.) 

In the simply laced cases, i.e., $\Phi$ of type $A_n$, $D_n$, $E_6$, $E_7$ or $E_8$, we 
scale $\Phi$ in such a way that for a  root $\alpha$  one has $(\alpha,\alpha)=2$ and we let
$y$ be a root. Then $(y,\alpha)\in\{0,\,\pm 1,\,\pm 2\}$ for $\alpha \in \Phi$. 
For $i\in \{0,\,\pm 1,\,\pm 2\}$ let $n_i=|\{\alpha \in \Phi\mid (y,\alpha)=i\}|$.
One has the obvious relations $n_{-i}=n_i$, $n_2=1$, and $n_0+2n_1+2n_2=|\Phi|$.
Furthermore (see~\cite{Bourbaki-Lie456}, Chap.~VI, \S1.11, Prop.~32), 
one has $|\Phi|=nh$ and $n_1=2h-4$, where
$h$ is the Coxeter number of~$\Phi$. Writing the polynomials $R_4$ and $R_6$ in the form
$R_l=\sum_{k=0}^{l/2}c_{l,k}(x,y)^{2k}(||x||/\sqrt{2})^{l-2k}$ one obtains
$$  \sum_{\alpha\in \Phi}R_l(\alpha(h^1),\,\ldots,\,\alpha(h^n))
=\sum_{i=-2}^2\sum_{k=0}^{l/2}n_ic_{l,k}i^{2k}=\cases{
                                 \frac{4 (h (n-10)+6 (n+2))}{n+2},  & for $l=4$, \cr
          \frac{4 \left(30 \left(n^2-16\right)+h \left(n^2-48 n+272\right)\right)}{n^2+12 n+32},   & for $l=6$.}$$ 
For $\Phi$ of type $A_n$, one has $h=n+1$ and therefore $R_4= \frac{4 \left(n^2-3 n+2\right)}{n+2}\not= 0$ 
for $n\not=1$, $2$ and $R_6=12\not=0$ for $n=2$. For the type $D_n$ one has $h=2n-2$ and 
thus $R_4= \frac{8 (n-4)^2}{n+2}\not=0$
for $n\not=4$ and  $R_6=24\not=0$ for $n=4$. For the types $E_6$ and $E_7$ one has  $h=12$ and $h=18$, 
which gives $R_4=0$, but $R_6=24\not=0$ and $R_6=\frac{192}{11}\not=0$, respectively.

Letting $\{e_1,\,\ldots,\,e_n\}$ be an orthonormal base, the remaining types of 
root systems~$\Phi$ can be realized as follows: $B_n$ by 
$\{\pm e_i,\,\pm e_i \pm e_j\mid 1 \leq i< j \leq n \}$ for $n\geq 2$,
$C_n$ by $\{\pm 2e_i,\,\pm e_i\pm e_j\mid 1 \leq i< j \leq n \}$ for $n\geq 3$,
$F_4$ by $\{\pm e_i,\,\frac{1}{2}(\pm e_1\pm e_2\pm e_3\pm e_4),\, \pm e_i\pm e_j\mid 1 \leq i< j \leq 4 \}$ 
and
$G_2$ by $\{\pm(e_1- e_2),\,\pm(e_2-e_3),\pm(e_3-e_1),\,\pm(2e_1-e_2-e_3),\,\pm(2e_2-e_1-e_3),\,\pm(2e_3-e_1-e_2)\}$.

For $\Phi$ of type $B_n$, we take $y=e_1+e_2$. By using the fact that the long roots of $\Phi$ form a root system
of type $D_n$ one obtains finally  $R_4= \frac{8 n^2-60 n+112}{n+2}\not=0$ for $n\not=4$ and $R_6=23\not=0$ for $n=4$.

For the type $C_n$, we take again $y=e_1+e_2$. Now the short roots of $\Phi$ form a root system
of type $D_n$ and one obtains finally  $R_4= \frac{8 \left(n^2-16\right)}{n+2}\not=0$ 
for $n\not=4$ and $R_6=-60\not=0$ for $n=4$.

For the type $F_4$, we take also $y=e_1+e_2$. Taking the additional $24$ short roots into account 
compared to $D_4$ one obtains $R_4=0$, but $R_6=21\not= 0$.

Finally, for the type $G_2$, we take $y=e_1-e_2$ and obtain $R_4=0$, but $R_6= -312 \not= 0$.
\phantom{xx}\eop

The last four Lemmas together and Theorem~\ref{highestweight} show that $V_1$ is
either abelian or simple of type $A_1$ or $E_8$.

\begin{prop}
The Lie algebra $V_1$ is not abelian.
\end{prop}
\pf Assume that $V_1$ is an abelian Lie algebra.
Let $V(\s)$ be the Heisenberg \voa of central charge $1$ generated by a one-dimensional
subalgebra $\s=\C\cdot  h$ of $V_1$, where $\langle h,\,h \rangle=2$. 
One checks directly that for 
$$v=(8h_{(-3)}h_{(-1)}-6 h_{(-2)}^2-2 h_{(-1)}^4){\bf 1}   \in V(\s)_4\subset V_4$$
one has $L_1'v=L_2'v=0$, where $L_n'=\frac{1}{2}{:}\sum_{k\in \Z} h_{(n-k)} h_{(k)}{:}$, $n\in {\bf Z}$,
are the generators for the Virasoro algebra of $V(\s)$. It follows as at the end of
the proof of Lemma~\ref{cartan-hw} that $v$ is also a Virasoro highest weight vector 
for the Virasoro algebra of $V$.

Using the associativity relation, one finds for the trace of $\o(v)$ on $V_1$:
\begin{eqnarray*}
\tr|_{V_1}\o(v) & = & 8\,\tr|_{V_1}\o(h_{(-3)}h_{(-1)}{\bf 1})-6\, \tr|_{V_1}\o(h_{(-2)}^2{\bf 1})
-2\, \tr|_{V_1}\o( h_{(-1)}^4{\bf 1}) \\
&=&8\,\tr|_{V_1}(3h_{(-1)}h_{(1)}+ h_{(0)}^2)-6\, \tr|_{V_1} h_{(0)}^2
-2\, \tr|_{V_1}( h_{(0)}^4+ 4 h_{(-1)}h_{(1)}h_{(0)}^2 )\\
& = & 24 \langle h,\,h \rangle \not=0.
\end{eqnarray*}

It follows from Theorem~\ref{highestweight} that $V_1$ cannot be a $6$-design.
\eop

We have proven the following result:

\begin{prop}\label{weylinv}
The Lie algebra $V_1$ is either  isomorphic to a Lie algebra of type 
$A_1$ or $E_8$. \eop
\end{prop}
 
\smallskip

\begin{prop}\label{m1charge}
The central charge of $V$ is either $1$, $8$ or $16$.
\end{prop}
\pf Let $\omega'$ be the
Virasoro element of a central charge $e=1$  Heisenberg \subvoa generated by a one-dimensional
subspace of the Cartan subalgebra~${\bf h}$ of~$V_1$.

Since $V_{1}$ is assumed to be a conformal \hbox{$6$-design},
one can apply Proposition~\ref{MODcond} to the 
$U_{\omega'}\otimes U_{\omega-\omega'}$-module  of conformal weight $h=1$
generated by $V_1$.
If $c\not=\frac{1}{2}$, $\frac{34}{35}$,~$1$, the proposition gives  $c=8$ or $c=16$ for the central charge.

\smallskip

The cases $c=\frac{1}{2}$ or $\frac{34}{35}$ can easily be excluded:

For $c=\frac{1}{2}$ one has $L_{1/2}(0)\subset V$. 
The only possible extension of $L_{1/2}(0)$ is by the simple current  $L_{1/2}(\frac{1}{2})$
of conformal weight $\frac{1}{2}$. In both cases one has $V_1=0$,  
a contradiction.

The central charge $c=\frac{34}{35}$ is excluded as in subsection~\ref{min0.5}.
\eop

The proof of Theorem~\ref{classfic7}, part (b) can now be finished.
As mentioned before, the \subvoa of $V$ generated by $V_1$ is an integrable
highest weight representation of the affine Kac-Moody algebra associated to $V_1$.
Such a representation is for simple $V_1$ determined by a level~$k$, where $k$ is a 
positive integer and the central charge of $\langle V_1\rangle $ and hence, 
by Lemma~\ref{samecharge}, of $V$  is given by $c=\frac{k\cdot \dim V_1}{k+\check g}$ 
where $\check g$ is the dual Coxeter number of $V_1$~\cite{FreZhu}. 
For $V_1$ of type $A_1$ one
has $c=\frac{3\,k}{k+2}$ and for $V_1$ of type $E_8$ one has $c=\frac{248\,k}{30+k}$.
By using Proposition~\ref{m1charge}, one finds that $k=1$ and $c=1$ for $V_1$ of type $A_1$ and
$k=1$ and $c=8$ for $V_1$ of type $E_8$ are the only possibilities.

For $k=1$, the \voa associated to the level~$1$ representation of an affine Kac-Moody algebra 
of type $A_1$ or $E_8$ is isomorphic to the lattice \voa $V_{A_1}$ or $V_{E_8}$ associated 
to the root lattice $A_1$ or $E_8$, respectively.

It follows that $V$ equals $V_{A_1}$ or $V_{E_8}$ since $V_{A_1}$ has besides $V_{A_1}$
only one irreducible module of conformal weight $1/4$ and $V_{E_8}$ is self-dual,
i.e., both \voas cannot be extended.
\smallskip

This finishes the proof of Theorem~\ref{classfic7}, part (b).


\subsection{Minimal weight $\frac{3}{2}$}\label{min1.5}  

Let $V$ be a \svoa as in part (c) of Theorem~\ref{classfic7} and let $V_{(0)}$ be
the even \subvoa. By assumption, $\dim V_2 \geq 2$, i.e., the minimal weight
of the  even \subvoa $V_{(0)}$ of $V$ is $2$.

Since we assumed that $V$ has a real form such that the invariant bilinear form on $V_2$
is positive definite,  it follows from~\cite{MeNe}, and~\cite{mia}, that the Virasoro
element $\omega\in V_2$ can be decomposed into the sum of two nonzero elements $\omega'$,
$\omega-\omega'\in V_2$  such that, after dividing by a factor~$2$, one has
two commuting idempotents of the algebra $V_2$.
The two elements generate commuting Virasoro vertex operator subalgebras $U_{\omega'}$ 
and $U_{\omega-\omega'}$  of $V$ of
central charge $e=2( \omega',\omega') $ and $c-e$, 
respectively, with  $0<e<c$ (cf.~\cite{Ho-dr}, Thm.~1.2.2).  

\smallskip

For $c\not=\frac{1}{2}$, $\frac{34}{35}$, $1$,
one can apply Proposition~\ref{MODcond} to the module $V_{(1)}$ of $V_{(0)}$
since $V_{3/2}$ is assumed to be a conformal $6$-design.
For $h=\frac{3}{2}$ one obtains $c=16$ or $c=23\frac{1}{2}$ for the central charge.

\medskip

The cases not covered by Proposition~\ref{MODcond} are again easily excluded:

For $c=\frac{1}{2}$ one has $L_{1/2}(0)\subset V $. 
There are no irreducible modules of $L_{1/2}(0)$ of conformal weight $\frac{3}{2}$.
Hence $V\cong L_{1/2}(0)$ and the minimal weight of $V$ cannot be $\frac{3}{2}$.

The central charge $c=\frac{34}{35}$ is excluded as in subsection~\ref{min0.5}.

If $c=1$, the central charges of the Virasoro algebras generated by
$\omega'$ and $\omega-\omega'$ both must be $\frac{1}{2}$, since $\frac{1}{2}$ is the 
smallest possible central charge in the unitary minimal series.
Thus $V$ contains a subalgebra isomorphic to $L_{1/2}(0)^{\otimes 2}$. 
There are no irreducible modules of $L_{1/2}(0)^{\otimes 2}$ of conformal weight $\frac{3}{2}$.
Hence $V\cong L_{1/2}(0)^{\otimes 2}$ and the minimal weight of $V$ would be $2>\frac{3}{2}$,
a contradiction.

\smallskip

This finishes the proof of Theorem~\ref{classfic7} part (c).


\subsection{Minimal weight $2$}\label{min2}  

Let $V$ be a \svoa as in part (d) of Theorem~\ref{classfic7}. All the
assumptions for $V$ also hold for $V_{(0)}$, in particular,
the minimal weights of $V$ and $V_{(0)}$ are the same.
If $V$ is not a \voa, we can replace therefore $V$ by $V_{(0)}$.

As the minimal weight of $V$ is $2$, the dimension $d=\dim  V_2$ of the degree-$2$-part
of $V$ is at least $2$.
Since we assumed that $V$ has a real form such that the invariant bilinear form on $V_2$
is positive-definite, we can find as in subsection~\ref{min1.5} an 
element $\omega'$ in $ V_2$ generating a Virasoro algebra of central charge $e$
such that $0<e<c$. 

\smallskip
Assume first that $c\not=\frac{1}{2}$, $\frac{34}{35}$, $(55\pm \sqrt{33})/2$, $1$.
In this case we can apply Proposition~\ref{m2VOAcond} and we have
$d=\frac{c(2388+955\,c+70\,c^2)}{2(748-55\,c+c^2)}$.

Since $V$ is assumed to be rational, the central charge $c$ has  
to be a rational number (see~\cite{DLM-twist}, Thm.~11.3).
\begin{lemma}\label{diophant}
The only positive rational numbers $c$ for which 
$$d=\frac{c(2388+955\,c+70\,c^2)}{2(748-55\,c+c^2)}$$ is a positive integer
are  $c= \frac{1}{2}$, $8$, $\frac{52}{5}$, $16$, $\frac{132}{7}$, $20$, $\frac{102}{5}$, 
$\frac{748}{35}$, $\frac{43}{2}$, $22$, $\frac{808}{35}$, $\frac{47}{2}$, $24$, 
$\frac{170}{7}$, $\frac{49}{2}$, $\frac{172}{7}$, $\frac{152}{5}$, $\frac{61}{2}$, 
$\frac{154}{5}$, $\frac{220}{7}$, $\frac{63}{2}$, $32$, $\frac{164}{5}$, $\frac{236}{7}$, $34$, 
$\frac{242}{7}$, $36$, $40$, $\frac{204}{5}$, $44$, $\frac{109}{2}$, $\frac{428}{7}$, $68$, 
$ \frac{484}{7}$, $\frac{187}{2}$, $132$, $1496$. 
\end{lemma}

\pf Let $c=\frac{p}{q}$ with coprime integers $p$ and $q$. The equation for $d$ can
be rewritten as
$$2\,q(748\,q^2-55\,pq+p^2)d=p(2388\,q^2+955\,pq+70\,p^2).$$
It follows successively $q|p(2388\,q^2+955\,pq+70\,p^2)$, $q|2388\,q^2+955\,pq+70\,p^2$,
$q|70\,p^2$, $q|70$, i.e., $c=k/70$ with a positive integer $k$. This gives
$$d=\frac{4805}{2} +  \frac{k}{2}+
 \frac{ -17611286000 + 15001210\, k}{2\,( 3665200 - 3850\,k + k^2 )}.$$
As $|( -17611286000 + 15001210\, k)/(2\,( 3665200 - 3850\,k + k^2 ))|< \frac{1}{2}$
for $k>70\, (107179 + \sqrt{11483743153}) \simeq 15\,003\,885.97\ldots$, the result follows 
by computing $d$ for all $k<15\,003\,886$. \eop

\smallskip

If the \voa $V$  is a self-dual \voa,
then $c$ has to be an integer divisible by $8$ 
(see~\cite{Ho-dr}, Cor.~2.1.3 and Section~\ref{extremal} above). 
From  Lemma~\ref{diophant} it follows that in this 
case $c=8$, $16$, $24$, $32$, $40$ or $1496$.
For self-dual \voas of central charge $8$ and $16$~\cite{Ho-dr}, Thm.~2.1.2 gives 
$\dim V_1=248$ and $\dim V_1=496$, respectively,
i.e., these two cases are impossible as the minimal weight of $V$ would be~$1$. 

If the \voa $V$ is not self-dual, then there exists an irreducible 
$V$-module $W$ different from $V$.
Proposition~\ref{MODcond} gives $h= \frac{124 + 31\,c \pm{\sqrt{31}}\,{\sqrt{368 + 24\,c - c^2}}}{496}$.
Since $V$ is assumed to be rational, the conformal weight $h$ of an irreducible module has 
to be a rational number (\cite{DLM-twist}, Thm.~11.3). A direct verification gives:
\begin{lemma}\label{rational}
The only values of $c$ listed in Lemma~\ref{diophant} for which 
\[ h= \frac{124 + 31\,c \pm{\sqrt{31}}\,{\sqrt{368 + 24\,c - c^2}}}{496} \] 
is rational are  $c= \frac{1}{2}$, $8$, $16$, $\frac{808}{35}$, $\frac{47}{2}$, $\frac{164}{5}$, 
$\frac{236}{7}$, $\frac{242}{7}$. \eop 
\end{lemma}

By applying Propositions~\ref{m2VOAcond} and~\ref{MODcond} and Lemmata~\ref{diophant} 
and~\ref{rational}, it follows therefore that, if $V$ is not self-dual, the only possible values for the 
central charge are $c=8$, $16$, $\frac{808}{35}$, $\frac{47}{2}$, $\frac{164}{5}$, $\frac{236}{7}$ or $\frac{242}{7}$.
One also obtains the values for $\dim V_2$ and $h$ as listed in Table~\ref{listem2}.

\medskip

The remaining cases for $c$ can again easily be excluded:

For $c=\frac{1}{2}$, one has $L_{1/2}(0)\subset V $. 
There are no irreducible modules of $L_{1/2}(0)$ of conformal weight $2$.
Hence $V\cong L_{1/2}(0)$ and the minimal weight of $V$ is larger than~$2$,
a contradiction.

For $c=\frac{34}{35}$, the same argument as in subsection~\ref{min0.5} holds.

If $c=1$, $V$ must contain a subalgebra isomorphic to $L_{1/2}(0)^{\otimes 2}$
by the argument given in subsection~\ref{min1.5}.
There are no irreducible modules of $L_{1/2}(0)^{\otimes 2}$ of conformal weight $2$.
Hence $V\cong L_{1/2}(0)^{\otimes 2}$. Let $\omega'$ be the Virasoro element of
 one subalgebra $L_{1/2}(0)$.
Let $d_h$ be the multiplicity of the eigenvalue $h$ of $L_0'=\omega_{(1)}'$ acting on $V_2$.
We have $d_0=1$ and  $d_{1/2}=d_{1/16}=0$. 
The vectors $v^{(2)}$ and $v^{(4)}$ given in~(\ref{v2}) and~(\ref{v4}) are
still well-defined Virasoro highest weight vectors. For the traces of $\o(v^{(4)})$ on $V_2$ one gets
$$ 
\tr|_{V_2} \o(v^{(4)})=  \frac{27}{5} + \frac{22}{5}\,{d_0} - 
    \frac{59}{5}\,{d_{1/2}} + \frac{571}{320}\,{d_{1/16}} $$
and so $\tr|_{V_2} \o(v^{(4)})=\frac{49}{5}$, a contradiction
to the conformal $6$-design property of $V_2$.

The cases $c=(55\pm \sqrt{33})/2$ are excluded because $c$ is not rational. 

\smallskip

This finishes the proof of Theorem~\ref{classfic7} part (d).


\subsection{Conformal $8$-designs}\label{8designs}

For a \voa with a module whose lowest degree part is a conformal $8$-design one has in
addition to Proposition~\ref{MODcond}:
\begin{prop}\label{MOD8cond}
Let $V$ be a \voa of central charge $c\not =-31$, $-\frac{44}{5}$, $-\frac{184}{105}$, 
$\frac{6}{55}$, $(-47 \pm 5\sqrt{57})/4$, $1$ and assume
there exists a module $W$ of $V$ of conformal weight~$h\not= 0$ 
such that the lowest degree part $W_h$ is a conformal $8$-design.
If there exist elements $\omega'$, $\omega-\omega' \in V_2$ 
generating two commuting Virasoro \voas of central charge $e$ and $c-e$, respectively, 
with $e$, $c-e \not = -\frac{46}{3}$, $-\frac{68}{7}$, $-\frac{22}{5}$,  $-\frac{3}{5}$, $0$, 
$\frac{1}{2}$ then
$$ (c-24 h+12) \bigl(10\,c^3 + (141-615 h)\,c^2 + 2(5740 h^2-3321 h+171)\,c \qquad\qquad\qquad\qquad $$
\begin{equation}\label{MOD8equation}
\qquad\qquad\qquad\qquad\qquad\qquad\qquad\qquad {}-24 (2870 h^3-2870 h^2+451 h+15 )\bigr)=0.
\end{equation}
If  $e=\frac{1}{2}$ or $e=c-\frac{1}{2}$ and one assumes that the Virasoro \voa generated by $\omega'$
or $\omega-\omega'$, respectively, is simple and the other assumptions hold,
then
$$ 152700 c^6+(3535420-2546820 h) c^5+2 (5519040 h^2-33944388 h+10007663) c^4 $$
$$ +(-66228480    h^3+505357184 h^2-303807330 h+24963561) c^3 \qquad\qquad  $$
$$ +(-2634772224  h^3+409756928 h^2+611923251 h-162937170) c^2 $$ 
$$ \qquad +3 (4450030592 h^3-4570094080 h^2+1282098891 h-49633193) c $$ 
\begin{equation}\label{MOD8equation2}
\qquad\qquad\qquad  +18 (120063488 h^3-120063488 h^2+34154597 h-1236817)=0. 
\end{equation}
\end{prop}
\pf The proof is similar to the proof of Proposition~\ref{MODcond}.

Assume first that $e$, $c-e\not=\frac{1}{2}$.   The conditions on the central charges
guarantee that the Virasoro \voas  generated by
 $\omega'$ and $\omega-\omega'$ are isomorphic to $M(x,0)/M(x,1)$ up to degree~$8$, 
where $x=e$ or $c-e$, respectively. Lemma~\ref{hwnumbers} shows that one can find
three linear independent highest weight vectors 
$v_a^{(8)}$, $v_b^{(8)}$ and  $v_c^{(8)}\in U_{\omega'}\otimes U_{\omega-\omega'}$.
By assumption, the degree~$h$ subspace $W_h$ is a conformal $8$-design. Thus one has 
the trace identities
\begin{equation}\label{modtraces8}
\tr|_{W_h}\o(v^{(2)})=\tr|_{W_h}\o(v^{(4)})=\tr|_{W_h}\o(v_a^{(8)})=\tr|_{W_h}\o(v_b^{(8)})=\tr|_{W_h}\o(v_c^{(8)})=0
\end{equation}
which form a homogeneous system of linear equations for $d^*=m_0^*$, $m_1^*$, $m_2^*$
$m_3^*$, and~$m_4^*$, where $d^*=\dim  W_h$ and $m_i^*=\tr|_{W_h} a_0^i$ for $i=1$, $2$, $3$, $4$.
Since $d^*>0$, the system has to be singular and the determinant has to vanish. 
This condition gives the proposition in this first case.

\smallskip
The case $e =\frac{1}{2}$ or $e=c-\frac{1}{2}$ is handled similar as in the proof
of Proposition~\ref{MODcond} by choosing an appropriate highest weight vector 
$v^{(8)}\in U_{\omega'}\otimes U_{\omega-\omega'}$ besides $v^{(2)}$ and $v^{(4)}$. \eop

\smallskip

Let now $V$ be a \svoa of minimal weight $m=1$ or $m=\frac{3}{2}$ satisfying the conditions
of Theorem~\ref{classfic8}~(i).
The proof of Theorem~\ref{classfic7} shows that Proposition~\ref{MOD8cond} is applicable
in these cases with $h=1$ or $h=\frac{3}{2}$, respectively. 
First, we assume  $e$, $c-e\not =\frac{1}{2}$.
For $h=1$, equation~(\ref{MOD8equation}) gives $c=12$ or 
$c=(177\pm \sqrt{22009})/10$; for $h=\frac{3}{2}$, one gets 
$c=24$, $10.04101...$, $19.13162...$, $48.97735...$ . These
central charges are impossible by Theorem~\ref{classfic7}, part (b) and (c).
Now we assume  $e =\frac{1}{2}$ or $e=c-\frac{1}{2}$.
For $h=1$, equation~(\ref{MOD8equation2}) has the positive real solutions $c=12$ and
$c=2.268296...$; for $h=\frac{3}{2}$, one gets the  positive real solutions
$c=3.342825...$ and $c=18.81561...$ . Again, these central charges are 
excluded by Theorem~\ref{classfic7}, part (b) and~(c).

\smallskip

Thus we have proven Theorem~\ref{classfic8}~(i) for $m=1$ or $m=\frac{3}{2}$.
The case $m=\frac{1}{2}$ follows from Theorem~\ref{classfic7}~(a).

\medskip

In addition to Proposition~\ref{m2VOAcond}, one has:
\begin{prop}\label{m2VOA8cond}
Let $V$ be a \voa of central charge $c\not= -31$, $-\frac{44}{5}$, $-\frac{184}{105}$,
 $\frac{6}{55}$, $36$, $\frac{47\pm 5\sqrt{57}}{4}$, $17.58127...$, $25.84832...$, $65.47039...$, $1$
with $V_1=0$ such that $V_2$ forms a conformal $8$-design.
If there exist elements $\omega'$, $\omega-\omega' \in V_2$ 
generating two commuting Virasoro \voas of central charge $e$ and $c-e$, respectively, 
with $e$, $c-e \not = 
-\frac{46}{3}$, $-\frac{68}{7}$, $-\frac{22}{5}$,  $-\frac{3}{5}$, $0$, $\frac{1}{2}$,
then 
$$d=  \frac{15 c \left(155 c^3+4133 c^2+32074 c+88392\right)}{20 c^3-2178
    c^2+65956 c-595056}.$$ 
If  $e=\frac{1}{2}$ or $e=c-\frac{1}{2}$ and one assumes that the Virasoro \voa generated by $\omega'$
or $\omega-\omega'$, respectively, is simple and also
$c\not=-26.45283...$, $0.23416...$,  $25.60637...$, $15.35332... \pm 13.56755... i$, 
and the other assumptions hold then
$d= -\frac{15 c \left(5734920 c^5+59136716 c^4+283246086 c^3+2858841411
    c^2+7908127017 c-2179288566\right)}{260520 c^5-7840184 c^4-72048858
    c^3+5528559692 c^2-75371626638 c+17347413996}.$
\smallskip
\end{prop}
\pf The proof is similar to the proof of Proposition~\ref{m2VOAcond}.

Assume first that $e$, $c-e\not=\frac{1}{2}$.   The conditions on the central charges
guarantee again that there
exist three linear independent highest weight vectors 
$v_a^{(8)}$, $v_b^{(8)}$ and  $v_c^{(8)}\in U_{\omega'}\otimes U_{\omega-\omega'}$.
By assumption, $V_2$ is a conformal $8$-design. Thus one has
the trace identities
$$
\tr|_{V_2}\o(v^{(2)})=\tr|_{V_2}\o(v^{(4)})=\tr|_{V_2}\o(v_a^{(8)})=\tr|_{V_2}\o(v_b^{(8)})=\tr|_{V_2}\o(v_c^{(8)})=0
$$
which form an inhomogeneous system of linear equations for $d=m_0$, $m_1$, $m_2$
$m_3$, and~$m_4$, where $d=\dim  V_2$ and $m_i=\tr|_{V_2} a_0^i$ for $i=1$, $2$, $3$, $4$.
The assumptions on the central charge guarantee that the system is non-singular
since its determinant is non-zero.
The solution for $d$ is the one given in the proposition and does not depend on $e$.

\smallskip
If $e =\frac{1}{2}$ or $e=c-\frac{1}{2}$, we
choose again an appropriate highest weight vector 
$v^{(8)}\in U_{\omega'}\otimes U_{\omega-\omega'}$ besides $v^{(2)}$ and $v^{(4)}$
and the result follows as in the proof of Proposition~\ref{m2VOAcond}.\eop

Let now $V$ be a \svoa of minimal weight $m=2$ satisfying the conditions
of Theorem~\ref{classfic8}~(i). 
As explained in subsection~\ref{min2},
Proposition~\ref{m2VOAcond} and Proposition~\ref{m2VOA8cond} are applicable
provided the conditions on the central charge of $V$ are satisfied. 

Assume first $e$, $c-e\not=\frac{1}{2}$. 
The two expressions given in these two propositions for $d$ together 
form an equation for $c$
with the solutions $c=-\frac{516}{13}$, $-\frac{44}{5}$ $-\frac{22}{5}$, $0$,
$24$, $\frac{142}{5}$. The case $c=\frac{142}{5}$ is impossible because 
$d=-164081<0$. Also, $c$ has to be positive.
The cases $c=\frac{1}{2}$, $\frac{34}{35}$ and $1$ are
excluded as in subsection~\ref{min2}. 
For $c=\frac{6}{55}$, $\frac{47\pm 5\sqrt{57}}{4}$, 
$17.58127...$, $25.84832...$, $65.47039...$, 
Proposition~\ref{m2VOAcond} gives a non-integer value for $d$.

This leaves the case $c=36$, which can be excluded by using the
non-singular linear system
$$\tr|_{V_2}\o(v^{(2)})=\tr|_{V_2}\o(v^{(4)})=\tr|_{V_2}\o(v_a^{(6)})=\tr|_{V_2}\o(v_a^{(8)})=\tr|_{V_2}\o(v_b^{(8)})=0$$
leading to $d=-67770<0$.

If $e=\frac{1}{2}$ or $e=c-\frac{1}{2}$, one has again two equations for $d$ and
one finds for $c$ the real solutions $c=-\frac{22}{5}$, $0$, $\frac{1}{2}$, $24$, $\frac{142}{5}$, 
$-8.45952...$ .
For the cases $c=0.23416...$,  $25.60637...$ with $c>0$ which are not yet excluded,
Proposition~\ref{m2VOAcond} gives non-integer values for $d$.

\smallskip

This finishes the proof Theorem~\ref{classfic8}~(i) for $m=2$ since for $c=24$ 
Proposition~\ref{m2VOAcond} gives $d=196884$.

\medskip

For part (ii) of Theorem~\ref{classfic8}, note that under the  extra assumptions there
Theorem~\ref{classfic7} (d) shows that $V$ has to be a self-dual
\voa of central charge~24. As discussed in section~\ref{extremal}, the character of $V$ 
has therefore to be equal to the character of the Moonshine module $V^{\natural}$.


\subsection{Known examples of conformal $6$-designs}\label{candidates}  

We will discuss which \svoas satisfying the conditions
of Theorem~\ref{classfic7} are known to support conformal $6$-designs.
We compare our results for \svoas with the analogous results for binary linear codes and 
integral lattices due to Lalaude-Labayle~\cite{LaLa-codesign} and Martinet~\cite{Ma-ladesign}.

\medskip

In the case of minimal weight $m=\frac{1}{2}$, the only example is
the self-dual \svoa $V_{\rm Fermi}$.
The homogeneous subspaces of 
$V_{\rm Fermi}$ are trivial conformal $t$-designs for all $t\geq 1$ because
$V_{(0)}$ is equal to the Virasoro highest weight module $L_{1/2}(0)$.

\smallskip

The only self-orthogonal binary code of minimal weight~$2$, whose set of minimal-weight
words supports a $3$-design 
is the trivial example of the self-dual code $C_2\cong \{(0,0),\,(1,1)\}$.

The only integral lattice of minimal norm $1$, whose set of minimal
vectors forms a spherical $7$-design 
is the one-dimensional lattice ${\bf Z}$ of integers.
In fact, the two vectors of any fixed positive integer length form a 
spherical $t$-design for all $t$.

\medskip

For $m=1$, the only examples are the  \voas $V_{A_1}$ and $V_{E_8}$.
As shown in Example~\ref{a1}, all the homogeneous subspaces of $V_{A_1}$ 
are conformal $t$-designs for arbitrary $t$.
As shown in Example~\ref{e8} and also in Example~\ref{extexa}, 
all the homogeneous subspaces of $V_{E_8}$ 
are conformal $7$-designs. However, Theorem~\ref{classfic8} (i) shows that $(V_{E_8})_1$
is not a conformal $8$-design.

\smallskip

The only self-orthogonal binary code of minimal weight~$4$, whose set of
minimal-weight words supports a $3$-design is the doubly-even self-dual Hamming 
code ${\cal H}_8$ of length~$8$. (One may also consider the trivial code $C_2'= \{(0,0)\}$.)

The only integral lattices of minimal norm $2$, whose set of minimal
vectors are spherical $7$-designs are the 
root lattices $A_1$ of rank~$1$ and $E_8$ of rank~$8$.

\medskip

In the case of minimal weight $m=\frac{3}{2}$, the shorter Moonshine module 
$V\!B^{\natural}\cong V\!B^{\natural}_{(0)}\oplus  V\!B^{\natural}_{(1)}$
is a \svoa of central charge $23\frac{1}{2}$ whose homogeneous subspaces
are conformal $7$-designs by Example~\ref{shortmoonshine} and~\ref{shortmoonshine2}. 
It follows from~\cite{Ho-babymoon} that $V\!B^{\natural}_8$
contains a Baby Monster invariant non-zero Virasoro highest weight vector and 
Theorem~\ref{classfic8} (i) shows that $V\!B^{\natural}_{3/2}$ is not a conformal
$8$-design.

As seen in Example~\ref{laminated}, the homogeneous subspaces of the central charge~$16$
\voa $V_{\Lambda_{16}}^+$ as well as of the module $K(\frac{3}{2})$ of
conformal weight $\frac{3}{2}$ are conformal $7$-designs. The module $K(\frac{3}{2})$ 
is the direct sum of (all) irreducible $V_{\Lambda_{16}}^+$-modules $K_{\mu}$ 
of conformal weight $\frac{3}{2}$ which are  simple currents of order~$2$. 
(This can be proven by using that $V_{\Lambda_{16}}^+$ 
has the structure of a framed \voa {}\cite{DGH-virs,LaYa-framedrep}.) 
It is not clear (and seems unlikely) that the
individual modules $K_{\mu}$ give rise to conformal $6$-designs. If this is not the case,
they cannot be used to extend $V_{\Lambda_{16}}^+$ to a \svoa of the required type.

One may therefore conjecture that the shorter Moonshine module $V\!B^{\natural}$ is the only
possible example for a \svoa as in Theorem~\ref{classfic7} part (c).

\smallskip 

The only self-orthogonal binary code of minimal weight~$6$, whose set of
minimal-weight words supports a $3$-design is the self-dual shorter Golay code
of length~$22$.

The only integral lattice of minimal norm $3$ whose set of minimal
vectors forms a spherical $7$-design
is the $23$-dimensional unimodular shorter Leech lattice $O_{23}$.

\medskip

For $m=2$, the known examples are the \voas 
$V_{\Lambda_8}^+$ and $V_{\Lambda_{16}}^+$ of central charge 
$8$ and $16$ of Example~\ref{laminated};
the even part $V\!B^{\natural}_{(0)}$ of the shorter Moonshine module
of central charge~$23\frac{1}{2}$ from Example~\ref{shortmoonshine} and~\ref{shortmoonshine2}; 
the self-dual Moonshine module $V^{\natural}$ of central charge~$24$ 
from Example~\ref{moonshine} and~\ref{extexa};
and the known extremal self-dual \voas of central charge~$32$ from Example~\ref{extexa}. 
In these examples, all homogeneous subspaces
are conformal $7$-designs; in the case of the Moonshine module, they are even
conformal $11$-designs.  In all examples besides $V^{\natural}$, Theorem~\ref{classfic8}
gives that the subspace $V_2$ is not a conformal $8$-design;
the Griess algebra $V_2^{\natural}$ is not a conformal $12$-design
(see~\cite{DoMa-highermoon}, Thm.~3 and the following discussion).

For the other values of $c$,
no examples are known. The homogeneous
subspaces of extremal \voas $V$  of central charge~$40$ are by 
Theorem~\ref{extremaldesign} conformal $3$-designs and one has $\dim V_2=20620$. 
All known examples of such \voas are ${\bf Z}_2$-orbifolds of lattice \voas, 
where the lattice is an extremal rank~$40$ lattice with minimal squared 
length $4$. Such \voas contain a Virasoro element $\omega'$ generating 
a Virasoro algebra of central charge~$1/2$. If $V_2$ is a
conformal $7$-design, one obtains however $d_0=\frac{441768}{37}$ which is impossible
(cf.~\cite{Ma-design}, Table~3.3).
Similarly, one can show that with the assumption $L_{1/2}(0)\subset V$ the only central 
charges~$c$ given in Lemma~\ref{diophant}
for which $d_0$, $d_{1/2}$ and $d_{1/16}$ are nonnegative integers
are $c=\frac{1}{2}$, $8$, $16$, $20$, $23\frac{1}{2}$, $24$, $24\frac{1}{2}$, $24\frac{4}{7}$,
$30\frac{2}{5}$, $30\frac{1}{2}$, $31\frac{1}{2}$, $32$, $32\frac{4}{5}$, $33\frac{5}{7}$,
and $36$.

\smallskip 

The only self-orthogonal binary codes of minimal weight~$8$, whose set of
minimal-weight words support a $3$-design are the doubly-even code
$C_8\cong \{(0,0,0,0,0,0,0,0)$, $(1,1,1,1,1,1,1,1)\}$ of length~$8$; 
the doubly-even Reed-Muller code ${\cal R}(1,4)$ of length~$16$; 
the doubly-even subcode of the shorter Golay code of length~$22$;
the doubly-even self-dual Golay code of length~$24$ and
the five extremal doubly-even self-dual codes of length $32$.

The only integral lattices of minimal norm $4$, whose set of minimal
vectors are spherical $7$-designs are the even laminated lattices
$\Lambda_{8}$ and $\Lambda_{16}$ of rank $8$ and $16$; the even sublattice
of the shorter Leech Lattice $O_{23}$ of rank~$23$;
the even unimodular Leech lattice $\Lambda_{24}$ of rank~$24$ 
and the even unimodular extremal lattices of rank~$32$.

\smallskip 

One may therefore conjecture that the examples mentioned above 
are already all examples of \svoas for Theorem~\ref{classfic7} 
part~(d).

\providecommand{\bysame}{\leavevmode\hbox to3em{\hrulefill}\thinspace}
\providecommand{\MR}{\relax\ifhmode\unskip\space\fi MR }
\providecommand{\MRhref}[2]{%
  \href{http://www.ams.org/mathscinet-getitem?mr=#1}{#2}
}
\providecommand{\href}[2]{#2}

\end{document}